\documentclass[11pt]{article}
\usepackage{amsmath, amssymb, amsthm}
\date{}
\title{Bisections of graphs}
\author{
Choongbum Lee
\thanks{Department of Mathematics,
UCLA, Los Angeles, CA 90095, email: choongbum.lee@gmail.com. Research supported in part by
a Samsung Scholarship.}
\and
Po-Shen Loh
\thanks{Department of Mathematical Sciences, Carnegie Mellon University,
Pittsburgh, PA 15213, e-mail: ploh@cmu.edu.  Research supported in part by an NSA Young Investigators Grant.}
\and
Benny Sudakov
\thanks{Department of Mathematics,
UCLA, Los Angeles, CA 90095, email: bsudakov@math.ucla.edu. 
Research supported in part by NSF grant DMS-1101185, by AFOSR MURI grant FA9550-10-1-0569 and by a USA-Israel BSF
grant.}
}

\theoremstyle{plain}
\newtheorem{theorem}{Theorem}[section]
\newtheorem*{theorem*}{Theorem}
\newtheorem{proposition}[theorem]{Proposition}
\newtheorem{lemma}[theorem]{Lemma}
\newtheorem*{lemma*}{Lemma}
\newtheorem{corollary}[theorem]{Corollary}

\newtheorem{definition}[theorem]{Definition}
\newtheorem*{problem*}{Problem}

\newcommand{\pr}[1]{\mathbb{P}\left(#1\right)}
\newcommand{\E}[1]{\mathbb{E}\left[#1\right]}
\newcommand{\Var}[1]{\text{Var}\left[#1\right]}

\begin{document}
\maketitle

\begin{abstract}
A bisection of a graph is a bipartition of its vertex set in which
the number of vertices in the two parts differ by at most 1, and
its size is the number of edges which go across the two parts. In
this paper, motivated by several questions and conjectures of
Bollob\'as and Scott, we study maximum bisections of
graphs. First, we extend the classical Edwards bound on
maximum cuts to bisections. A simple corollary of
our result implies that every graph on $n$ vertices and
$m$ edges with no isolated vertices, and
maximum degree at most $n/3 + 1$, admits a bisection of size at least
$m/2 + n/6$. Then using the tools that we developed
to extend Edwards's bound, we prove a judicious bisection result
which states that graphs with large minimum degree have a bisection
in which both parts span relatively few edges. A special case
of this general theorem answers a conjecture of Bollob\'as and Scott,
and shows that every graph on $n$ vertices and $m$ edges of minimum
degree at least 2 admits a bisection in which the number of edges in
each part is at most $(1/3+o(1))m$. We also present
several other results on bisections of graphs.
\end{abstract}

%  Completely answer problem of Bollob\'as and Scott.
%  Challenge: a wide family of asymptotically extremal examples.
%  Key ingredient is a novel new concept of tight components.

\section{Introduction}

Many classical partitioning problems in Combinatorics and Computer
Science seek a partition of a combinatorial object (e.g., a graph,
directed graph, hypergraph, etc.) which optimizes a single parameter.
Two famous examples are Max Cut and Max Bisection, which have been widely
studied in Combinatorics and Computer Science (see, e.g., \cite{FrJe97, GoWi95,
  Hastad01, TSSW00} and their references). In these problems, the objective
is to find a bipartition which maximizes the number of crossing edges,
with Max Bisection having the additional constraint that the two
vertex classes should differ in size by at most one.

The most basic result in that area is that every graph with $m$ edges
has a cut of size at least $\frac{m}{2}$.  This
can be seen by considering the expected number of crossing edges in a
uniformly random partition, or by analyzing the following natural greedy
algorithm.  In this algorithm, one processes the vertices in an arbitrary
order, adding each new vertex to the side which maximizes the number
of edges crossing back to the previously positioned vertices. Moreover,
one can similarly obtain the same estimate for maximum bisection.

There have been many improvements of the above bound, most notable among
which is the following result of Edwards \cite{Edwards73} which is
tight, e.g., for complete graphs.  (Here and in the remainder, $e(U,
W)$ denotes the number of edges between two disjoint sets $U$ and
$W$.)
\begin{equation}
  \label{ineq:edwards}
  \max_{V_1 \cap V_2 = \emptyset} e(V_1, V_2)
  \geq
  \left\lceil
    \frac{m}{2} + \sqrt{\frac{m}{8} + \frac{1}{64}} - \frac{1}{8}
  \right\rceil \,.
\end{equation}
The area of \emph{judicious partitioning}\/ considers the general
family of problems in which one seeks to optimize multiple parameters
simultaneously.  In the Max Cut setting, the canonical example is the
beautiful result of Bollob\'as and Scott which shows that one can find
a cut $V_1 \cup V_2$ which not only achieves the Edwards bound
\eqref{ineq:edwards}, but also induces few edges in each $V_i$.
\begin{theorem}
  \label{thm:bs-judicious-cut}
  (Bollob\'as and Scott, see \cite{BoSc99}.)  Every graph $G = (V, E)$
  with $m$ edges has a bipartition $V = V_1 \cup V_2$ which achieves
  the bound in \eqref{ineq:edwards}, and each vertex class spans at
  most
  \begin{displaymath}
    \frac{m}{4} + \sqrt{ \frac{m}{32} + \frac{1}{256} } - \frac{1}{16}
  \end{displaymath}
  edges.
\end{theorem}
Note that this simultaneously optimizes three parameters: the number
of edges crossing the cut, and the number of edges inside each $V_i$ for $i=1,2$.
It is impossible to obtain results of this nature through
probabilistic arguments which are based on first-moment (expected
value) calculations alone.  Indeed, the most basic application of the
probabilistic method constructs a random variable $X$ to measure a
single parameter, and calculates its expectation for a random
construction. Then, as we did above for maximum cut, one can conclude that
there is an outcome in which $X$ takes a value which is at least its
expectation.

Part of the appeal of judicious partitioning problems thus stems from
the fact that they push the envelope for probabilistic methods: when
one needs to satisfy two separate parameters $X$ and $Y$
simultaneously, the above argument breaks down, because $X \geq \E{X}$
and $Y \geq \E{Y}$ do not necessarily hold at the same time.  It is
therefore no coincidence that many arguments in the area of judicious
partitioning blend techniques from probabilistic and extremal
combinatorics, using variance calculations or martingale concentration
inequalities to strengthen conclusions drawn from the initial
first-moment calculations.  This area has attracted many researchers,
who have now produced a wealth of results, covering combinatorial
structures spanning graphs, directed graphs, and hypergraphs.  We
direct the interested reader to any of \cite{BoReTh93, BoSc00, BoSc04,
  BoSc10, Ha11, KuOs07, MaYaYu10, MaYu11} (by no means a comprehensive
list), or to either of the surveys \cite{BoSc02a, Scott06}.

\subsection{Large bisections}

In this paper, we focus on bisections. Compared to cuts, bisections
are much more complicated to analyze.  For example, no bisection
analogue to Theorem \ref{thm:bs-judicious-cut} was previously known.
Also, in Computer Science, even long after Goemans and Williamson
\cite{GoWi95} discovered a 0.878-factor approximation to Max Cut via
semi-definite programming, the best Max Bisection results still lag
behind, from the 0.651-factor approximation of Frieze and Jerrum
\cite{FrJe97} that it inspired, to the current best 0.703-factor
approximation of Feige and Langberg \cite{FeLa06}.

Note that when $m$ is linear in $n$, the Edwards bound
\eqref{ineq:edwards} for cuts is relatively weak, giving only
\begin{displaymath}
  %\max_{V_1 \cap V_2 = \emptyset} e(V_1, V_2)
  \text{Max Cut}
  \geq
  \frac{m}{2} + \Omega( \sqrt{n} ) \,.
\end{displaymath}
As cuts have been more successfully studied than bisections, there
is an existing result which still has useful implications for sparse
graphs. As noticed in \cite{PoTu82, ErGyKo97},
the results that Edwards proved in \cite{Edwards75} implicitly
imply that connected graphs admit a bipartition of size at least
\begin{equation}
  \label{ineq:cut-connected}
  \text{Max Cut}
  \geq
  \frac{m}{2} + \frac{n-1}{4} \,.
\end{equation}
As a similar result, Erd\H{o}s, Gy\'arf\'as, and Kohayakawa \cite{ErGyKo97} established the
following bound for graphs without isolated vertices:
\begin{equation}
  \label{ineq:cut+n/6}
  \text{Max Cut}
  \geq
  \frac{m}{2} + \frac{n}{6} \,,
\end{equation}
which can be proved by using the probabilistic method to produce a suitable ordering of
the vertices to feed into the greedy algorithm alluded to in the
beginning of the Introduction.

Unfortunately, neither result transfers directly, because every
bisection of the star $K_{1,n-1}$ has size $\frac{n}{2} =
\frac{m+1}{2}$.  This situation can be circumvented by imposing
some minor conditions. To state our result, we introduce the
following new notion, with respect to which we can prove sharp
bounds on Max Bisection.

\begin{definition}
  \label{def:tight}
  Let $T$ be a connected component of a given graph.  We say that $T$
  is a \textbf{tight component} if it has the following properties.
  \begin{description}
  \item[(i)] For every vertex $v \in T$, the subgraph induced by $T
    \setminus \{v\}$ contains a perfect matching.
  \item[(ii)] For every vertex $v \in T$ and every perfect matching of
    $T \setminus \{v\}$, no edge of the matching has exactly one
    endpoint adjacent to $v$.
  \end{description}
\end{definition}

We can now define the parameter $\tau$ to be the number of tight
components in a graph.  Our first result is a very precise extension
of \eqref{ineq:cut-connected} and \eqref{ineq:cut+n/6} to bisections.

\begin{theorem}
  \label{thm:tight}
  Let $G$ be a graph of maximum degree $\Delta$ which has $\tau$ tight
  components.  Then there is a bisection of size at least
  \begin{displaymath}
    \frac{m}{2} + \frac{n - \max\{\tau, \Delta - 1\}}{4}.
  \end{displaymath}
\end{theorem}

This theorem is tight in both parameters $\tau$ and $\Delta$. For example
if the graph $G$ consists of $n/3$ vertex-disjoint copies of a triangle,
then $m= n, \tau = n/3, \Delta=2$, and one can easily see that the largest bisection
is of size $2n/3 = m/2 + n/6 = m/2 + (n - \tau)/4$.
As another example, if the graph is $K_{1,n-1}$ for an even integer $n$,
then $m = n-1, \tau = 0, \Delta = n-1$, and the largest bisection
is of size $n/2 = m/2 + 1/2 = m/2 + (n - \Delta + 1)/4$. Many other
tight examples can be constructed.

Analogues to \eqref{ineq:cut-connected} and \eqref{ineq:cut+n/6}
follow as easy corollaries from Theorem \ref{thm:tight}.

\begin{corollary}
  \label{cor:bisectionconnected}
  For $n \ge 2$, every $n$-vertex connected graph with $m$ edges and maximum degree
  $\Delta$ has a bisection of size at least $\frac{m}{2} + \frac{n + 1
    - \Delta}{4}$.
\end{corollary}

\begin{corollary}
  \label{cor:bisection+n/6-maxdegree}
  Every $n$-vertex graph with $m$ edges, maximum degree at most
  $\frac{n}{3} + 1$, and no isolated vertices, has a bisection of size
  at least $\frac{m}{2} + \frac{n}{6}$.
\end{corollary}

More importantly, the ideas used to prove Theorem \ref{thm:tight}
appear to be fairly robust.  Variants on the theme produce several
additional results, one of which provides a parameterized
interpolation between cut and bisection results.  Let an
\emph{$\alpha$-bisection}\/ be a cut in which both sides have at
least $(\frac{1}{2} - \alpha) n$ vertices.

\begin{theorem}
  \label{thm:almost-bisection}
  Let $0 \leq \alpha \leq \frac{1}{6}$ be a fixed
  parameter.  Then every $n$-vertex graph with $m$ edges and no isolated
  vertices contains an $\alpha$-bisection of size at least
  $\frac{m}{2} + \alpha n$.
\end{theorem}

\noindent \textbf{Remark.}\, This provides further insight to the
bound $\frac{m}{2} + \frac{n}{6}$ of \eqref{ineq:cut+n/6}, since it
shows that the bound can in fact always be achieved by a
$\frac{1}{6}$-bisection.

\subsection{Judicious bisections}

A \emph{judicious bisection} of a graph is a bisection in
which both parts span few edges. This is a qualitative definition,
and a quantitative
definition will be given for each problem that we consider.
There is a nontrivial obstacle to extending results from judicious
cuts to judicious bisections, which is the fact that the
direct analogue is false.  Indeed, if one considers the star
$K_{1,n-1}$, one observes that although the natural bipartition
produces an excellent judicious cut which spans no edges on either
side, every \emph{bisection}\/ spans exactly $\frac{n}{2} - 1$ edges
on one side, and none on the other.  This is far from the ideal of
$\big( \frac{1}{4} + o(1) \big) m$ edges which is suggested by an
expected value calculation for a random bisection.

To circumvent this issue, Bollob\'as and Scott proposed several
conditions under which they suspected the ideal result would hold.
Specifically, in their survey \cite{BoSc02a}, they suggested that
perhaps one of $\delta \rightarrow \infty$ or $\Delta = o(n)$ might be
enough.  (Here and in the remainder, $\delta$ and $\Delta$ correspond
to the minimum and maximum degree of the graph.)  It is clear that
these conditions exclude the star $K_{1,n-1}$, together with similar
graphs which also are obstructions.  Our first result proves their
conjecture.

%\begin{theorem}
%  \label{thm:bs-conjecture}
%  For any $\varepsilon > 0$, there exist constants $C, \gamma > 0$ such
%  that the following holds for sufficiently large $m, n$.  Let $G =
%  (V, E)$ be an $n$-vertex graph with $m$ edges, where either
%  \textbf{(i)} $m \geq Cn$, or \textbf{(ii)} all degrees are at most
%  $\gamma n$.  Then $G$ admits a bisection $V = V_1 \cup V_2$ in which
%  each $V_i$ spans at most $\big( \frac{1}{4} + \varepsilon \big) m$
%  edges.
%\end{theorem}
\begin{theorem}
  \label{thm:bs-conjecture}
  Let $\varepsilon$ be a fixed positive constant and let $G$ be an $n$-vertex
  graph with $m$ edges such that \textbf{(i)} $m \ge \varepsilon^{-2}n$, or
  \textbf{(ii)} all degrees are at most $\frac{\varepsilon^2 n}{2}$.
  If $n$ is sufficiently large, then $G$ admits a bisection
  $V = V_1 \cup V_2$ in which each $V_i$ spans at most
  $\big( \frac{1}{4} + \varepsilon \big) m$ edges.
\end{theorem}

Continuing the study of judicious bisection of graphs with large
minimum degree, Bollob\'as and Scott suggested that once the minimum degree was
at least 2, one could substantially improve upon the behavior of
the now-excluded pathological $K_{1, n-1}$ example.  Specifically,
they conjectured in \cite{BoSc02a} that such graphs always have
bisections spanning at most $\frac{m}{3}$ edges in each part.  This
result would be tight, for example in the case when the graph was a
triangle.  (See the discussion after Theorem \ref{thm:bisection-min}
for further examples.)  Our next theorem asymptotically confirms
Bollob\'as and Scott's conjecture.  Here and in the remainder, $e(U)$
denotes the number of edges spanned by the set $U$.

\begin{theorem}
  \label{thm:bisection-mindegreetwo}
  Every $n$-vertex graph $G = (V, E)$ with $m$ edges and minimum
  degree at least 2 has a bisection in which both
  \begin{displaymath}
    e(V_i) \leq \left( \frac{1}{3} + o(1) \right) m \,.
  \end{displaymath}
\end{theorem}

In \cite{BoSc02a, Scott06}, they further asked how the bound changes
as the minimum degree condition grows, but did not conjecture the
value of the constant.  By Theorem \ref{thm:bs-conjecture}, the
asymptotic fraction should shrink to $\frac{1}{4}$ with larger
$\delta$.  Our next result determines the asymptotic constants
precisely, providing a complete answer to the question of Bollob\'as and Scott.

\begin{theorem}
  \label{thm:bisection-min}
  Let $\delta$ be a positive even integer.  Every $n$-vertex graph $G
  = (V, E)$ with $m$ edges and minimum degree at least $\delta$ has a
  bisection in which both
  \begin{displaymath}
    e(V_i) \leq \left( \frac{\delta + 2}{4 (\delta + 1) } + o(1) \right) m \,.
  \end{displaymath}
\end{theorem}

\noindent \textbf{Remark.}\, For odd $\delta$, one can deduce a bound
by applying the above result for the even integer $\delta - 1$, and
this is asymptotically tight (see below).  The dependence between the
constant and $\delta$ is therefore rather non-trivial.

\vspace{2mm}

Adding to the challenge, it turns out that the extremal examples come
from at least two strikingly different families.  For the first
family, take the vertex-disjoint union of an arbitrary number of
cliques of order $\delta$, together with an arbitrary number of
cliques of order $\delta + 1$, and add a new vertex adjacent to all
vertices in the cliques.  For the second family, consider the complete
bipartite graph $K_{\delta+1, n - \delta - 1}$.  Note that we obtain
graphs with minimum degree $\delta + 1$ (odd) from the second family
of constructions, as well as from the first family of constructions
when only $K_{\delta+1}$ are used.  These then provide examples for
asymptotic tightness in the remark above.

\begin{proposition}
  \label{prop:bisection-min-tight}
  Constructions from both of the above families only admit bisections
  which span at least $\big( \frac{\delta+2}{4(\delta+1)} - o(1)
  \big)m$ edges on some side.
\end{proposition}

\vspace{2mm}

The main tool in proving Theorems~\ref{thm:bisection-mindegreetwo}
and \ref{thm:bisection-min} is the following strengthening of
Theorem~\ref{thm:bs-conjecture} for sparse graphs,
which can also be viewed as a ``randomized'' version of Theorem \ref{thm:tight}.
It builds upon the proof of Theorem \ref{thm:tight} and uses martingale concentration
techniques to produce a bisection which is not only large, but also
judicious.  The key benefit is its parameterization in terms of the
number of tight components.

\begin{theorem}
  \label{thm:tightrandom}
  Given any real constants $C, \varepsilon > 0$, there exist $\gamma, n_0
  > 0$ for which the following holds.  Every graph $G = (V, E)$ with
  $n \geq n_0$ vertices, $m \leq Cn$ edges, maximum degree at most
  $\gamma n$, and $\tau$ tight components, admits a bisection $V = V_1
  \cup V_2$ where each $V_i$ induces at most
  \begin{displaymath}
    \frac{m}{4} - \frac{n - \tau}{8} + \varepsilon n
  \end{displaymath}
  edges.
\end{theorem}

To make use of this theorem for graphs whose maximum degree is not
bounded, we first identify the large-degree vertices, say, of degree
at least $n^{3/4}$, and find an optimal split of these vertices.
Then armed with the appropriate parameterization, the remainder of
the proof follows from a (somewhat intricate) analysis of the
relationship between tight components, the number of edges,
large-degree vertices, and the minimum degree condition.

\subsection{Related bisection results}
\label{sec:related-results}

In the process of extending classical results from cuts to bisections, we also
discovered several related statements which we present
in this section.  The judicious bipartition bound in Theorem
\ref{thm:bs-judicious-cut} is tight for complete graphs, and therefore
cannot be improved in general.  This is clear: if the maximum cut only
has size $(\frac{1}{2} + o(1))m$, then no bipartition can induce
substantially fewer than $\frac{m}{4}$ edges on each side.  There is
only potential for improvement if the graph already has a cut which is
larger than the guaranteed $\frac{m}{2}$.  In this setting, Alon,
Bollob\'as, Krivelevich, and Sudakov \cite{AlBoKrSu03} proved that
when the maximum cut is large, then indeed there is a judicious
bipartition in which each side induces substantially fewer than
$\frac{m}{4}$ edges.

Bollob\'as and Scott \cite{BoSc10} asked whether the same is also
true for judicious bisections; we show that this is not the case.

\begin{proposition}
  \label{prop:big-bisection-no-judicious}
  For arbitrarily large $m$, there are $m$-edge graphs with bisections
  of size at least $\frac{3}{4} m$, but with $\max\{e(V_1), e(V_2)\}
  \geq \frac{m}{4}$ for every bisection $V_1 \cup V_2$.
\end{proposition}

Several researchers \cite{AnGrLi83, LeTu82, Locke82} also noticed
that if a graph has chromatic number at most $k$, then it admits a
bipartition of size at least $\frac{k+1}{2k} m$ when $k$ is odd and
$\frac{k}{2(k-1)} m$ when $k$ is even.  This unfortunately is
false for bisections, as the star $K_{1,{n-1}}$ has chromatic number
2, but maximum bisection $\frac{n}{2}$, not $\frac{2n}{3}$. Instead,
we prove an almost-bisection result under a maximum degree
condition.

\begin{theorem}
  \label{thm:max-degree-bisection}
  For any positive integer $r$, if $G$ is a graph with maximum degree
  at most $r$, then it has a bipartition of size at least
  $\frac{r+1}{2r} m$ when $r$ is odd and $\frac{r+2}{2(r+1)}m$ when $r$
  is even, in which the part sizes differ by at most $\frac{r}{2} + 1$.
\end{theorem}

An obvious corollary of Theorem \ref{thm:max-degree-bisection} is
that if $G$ is a graph with maximum degree at most $r$,
then it has a true bisection of size at least $\frac{r+1}{2r} m -
\frac{r(r+1)}{4}$ when $r$ is odd and $\frac{r+2}{2(r+1)} m -
\frac{r(r+2)}{4}$ when $r$ is even. We also give a simple proof of the
following result of Bollob\'as and Scott \cite{BoSc04}.

\begin{theorem}
  \label{thm:regular-bisection}
  For any positive integer $r$, every $r$-regular graph has a
  bisection of size at least $\frac{r+1}{2r} m$ when $r$ is odd and
  $\frac{r+2}{2(r+1)} m$ when $r$ is even.
\end{theorem}

\subsection{Notation and organization}

Since the majority of our results are asymptotic in nature, we will
implicitly ignore rounding effects (e.g., when bisecting graphs of odd
order) whenever these effects are of smaller order than our error
terms.  The following (standard) asymptotic notation will be utilized
extensively.  For two functions $f(n)$ and $g(n)$, we write $f(n) =
o(g(n))$, $g(n) = \omega(f(n))$ if $\lim_{n \rightarrow \infty}
f(n)/g(n) = 0$, and $f(n) = O(g(n))$ or $g(n) = \Omega(f(n))$ if there
exists a constant $M$ such that $|f(n)| \leq M|g(n)|$ for all
sufficiently large $n$.  We also write $f(n) = \Theta(g(n))$ if both
$f(n) = O(g(n))$ and $f(n) = \Omega(g(n))$ are satisfied.

This paper is organized as follows.  In the next section, we show
how Theorem \ref{thm:bs-conjecture}
follows from a rather short second-moment argument.  In Section
\ref{sec:large-bisection}, we study tight components, and develop
the main tools needed to prove
Theorems \ref{thm:tight} and \ref{thm:almost-bisection}, and
Corollaries \ref{cor:bisection+n/6-maxdegree} and
\ref{cor:bisectionconnected}. In Section \ref{sec:judicious-bisection},
we further combine these tools with martingale concentration
results and prove Theorem \ref{thm:tightrandom}.
We move to our main result in Sections
\ref{sec:bisection-min-2} and \ref{sec:bisection-min}, where we
first prove the $\delta = 2$ case (Theorem \ref{thm:bisection-mindegreetwo}), 
and then the general case
(Theorem \ref{thm:bisection-min}).  We further confirm
its asymptotic optimality by proving
Proposition \ref{prop:bisection-min-tight} in Section
\ref{sec:bisection-min-tight}.  Section \ref{sec:related-proofs}
contains the proofs of the various related results introduced in
Section \ref{sec:related-results}.  Finally, we conclude with some
remarks.
%and open problems.

\section{The basic second-moment argument}
\label{sec:bs-conjecture}

We begin by answering two questions of Bollob\'as and Scott, regarding
settings in which it was suspected that a judicious bisection would
exist.  These were the settings where (i) the minimum degree grew with
$n$, or (ii) the maximum degree was of order $o(n)$.  Either condition
would rule out the basic pathological counterexample, which was based
on a star.

It turns out that a second-moment argument suffices to answer these
questions.  The most elementary bisection algorithm constructs an
arbitrary pairing of the vertices, and then separates each pair across
$V_1 \cup V_2$, independently and uniformly at random.  We will now
show that this performs well in both of the above settings.

\begin{lemma}
  \label{lem:bisection-variance}
  Let $G = (V, E)$ be an $n$-vertex graph with $m$ edges, and define
  \begin{equation}
    \label{eq:var-lambda}
    \Lambda = \frac{1}{16} \left( 3m + \sum_{v \in V} d(v)^2 \right) \,,
  \end{equation}
  where $d(v)$ denotes the degree of vertex $v$.  Then $G$ admits a
  bisection $V = V_1 \cup V_2$ for which each
  \begin{displaymath}
    e(V_i) \leq \frac{m}{4} + \sqrt{2 \Lambda} \,.
  \end{displaymath}
\end{lemma}

\noindent \textbf{Proof.}\, Define the random variables $Y_1$ and
$Y_2$ to be the numbers of edges induced by each $V_i$ after the
elementary random bisection algorithm, described above.  It is clear that $\E{Y_1}
\leq \frac{m}{4}$, because for each edge $e = \{v, w\}$, $\pr{v, w
\in V_1} = 0$ if $v, w$ were initially paired, and it is
$\frac{1}{4}$ otherwise.

To calculate the variance of $Y_1$, define an indicator variable $I_e$
for each edge $e$, which takes the value 1 precisely when both
endpoints of $e$ fall in $V_1$.  Note that for any $e$ whose endpoints
were initially paired, $I_e$ is always 0.  Thus, such $e$ will never
contribute to $Y_1$; let $F$ be the set of all other edges (i.e.,
edges whose endpoints were not initially paired).  We then have
\begin{align*}
  \E{ \left(
      Y_1 - \E{Y_1}
    \right)^2 }
  &=
  \E{ \left(
      \sum_{e \in F} \left( I_e - \E{I_e} \right)
    \right)^2 } \\
  &=
  \E{ \sum_{e \in F} \left( I_e - \E{I_e} \right)^2
    +
    \sum_{e \neq f \in F} \left( I_e - \E{I_e} \right) \left( I_f - \E{I_f} \right) } \\
  &=
  \sum_{e \in F} \left( \E{I_e^2} - \E{I_e}^2 \right)
  +
  \sum_{e \neq f \in F} \left( \E{I_e I_f} - \E{I_e} \E{I_f} \right)
\end{align*}
Each term in the first sum is precisely $\frac{1}{4} - \big(
\frac{1}{4} \big)^2 = \frac{3}{16}$.  To estimate the second sum,
observe that $I_e$ and $I_f$ are independent whenever the two pairs
which hold the endpoints of $e$ are disjoint from the two pairs which
hold the endpoints of $f$.  In those cases, $\E{I_e I_f} - \E{I_e}
\E{I_f} = 0$.  Otherwise, the endpoints of $e$ and $f$ are contained
in either 2 or 3 pairs.  If only two pairs are involved, it is easy to
see that it is impossible for both of $e \neq f$ to lie in $V_1$, so
$\E{I_e I_f} - \E{I_e} \E{I_f} = -\big(\frac{1}{4}\big)^2$.  If three
pairs are involved but $e$ and $f$ do not share a vertex, it is also
impossible for both $e$ and $f$ to lie in $V_1$, so again $\E{I_e I_f}
- \E{I_e} \E{I_f} = -\big(\frac{1}{4}\big)^2$.  The final possibility
has $e$ and $f$ sharing a common endpoint, and their 3 total vertices
lie in 3 different pairs.  Then, the probability of $e$ and $f$ both
falling in $V_1$ is $\frac{1}{8}$, and hence $\E{I_e I_f} - \E{I_e}
\E{I_f} \leq \frac{1}{8} - \big( \frac{1}{4} \big)^2 = \frac{1}{16}$.
Therefore,
\begin{displaymath}
  \Var{Y_1}
  \leq
  \sum_{e \in F} \frac{3}{16}
  +
  \sum_{e \neq f \text{, incident}} \frac{1}{16}
  <
  \frac{3}{16} |F|
  +
  \frac{1}{16} \sum_{v \in V} d(v)^2
  \leq \Lambda \,,
\end{displaymath}
where the strict inequality comes from writing $d(v)^2$ instead of
$d(v) (d(v) - 1)$.  Chebyshev's inequality then gives
\begin{displaymath}
  \pr{ Y_1 \geq \E{Y_1} + \sqrt{2 \Lambda} }
  \leq
  \frac{\Var{Y_1}}{2 \Lambda}
  <
  \frac{1}{2} \,.
\end{displaymath}
The same inequality holds for $Y_2$ by symmetry, so we conclude that
there is an outcome in which both $e(V_i) \leq \frac{m}{4} +
\sqrt{2\Lambda}$.  \hfill $\Box$

\vspace{2mm}

We now proceed to prove Theorem \ref{thm:bs-conjecture}, which
considers the situations when (i) the number of edges is $\omega(n)$,
or (ii) the maximum degree is $o(n)$.

\vspace{2mm}

\noindent \textbf{Proof of Theorem \ref{thm:bs-conjecture}.}\, Let
$\varepsilon > 0$ be given.  We begin with (i).  We will show that
whenever there are $m \geq \varepsilon^{-2} n$ edges and $n$ is
sufficiently large, there is a judicious bisection in which each
$e(V_i) \leq \big( \frac{1}{4} + \varepsilon \big) m$.  Indeed, then the
trivial bound $d(v) < n$ applied to \eqref{eq:var-lambda} gives
\begin{displaymath}
  \Lambda
  <
  \frac{1}{16} \left( 3m + n \sum_{v \in V} d(v) \right)
  =
  \frac{1}{16} \left( 3m + n (2m) \right) \,,
\end{displaymath}
which is well below $\frac{1}{2} mn$ for large $n$.  Using $m \geq
\varepsilon^{-2} n$, we find that $\Lambda < \frac{1}{2} \varepsilon^2 m^2$,
which together with Lemma \ref{lem:bisection-variance} implies that
there is a bisection in which each $e(V_i) \leq \frac{m}{4} +
\sqrt{2 \cdot \frac{1}{2} \varepsilon^2 m^2}$, as desired.

For (ii), we instead assume that $d(v) \leq \frac{\varepsilon^2 n}{2}$.
%Although we do not assume that $m = \omega(n)$, we still assume that $n$ is sufficiently large.
Substituting our new bound for $d(v)$
into \eqref{eq:var-lambda}, we find that
\begin{displaymath}
  \Lambda
  <
  \frac{1}{16} \left( 3m + \left( \frac{\varepsilon^2 n}{2} \right) \sum_{v \in V} d(v) \right)
  =
  \frac{1}{16} \left( 3m + \varepsilon^2 n m \right) \,,
\end{displaymath}
which is well below $\frac{1}{8} \varepsilon^2 n m$ for large $n$.  If $m
\geq \frac{n}{4}$, then we have $\Lambda < \frac{1}{2} \varepsilon^2
m^2$, which completes the proof as for (i) above.  On the other hand,
if $m < \frac{n}{4}$, then there are at least $\frac{n}{2}$ isolated
vertices.  We may then apply the original result of Bollob\'as and
Scott (Theorem \ref{thm:bs-judicious-cut} in the Introduction) to the
graph induced by the non-isolated vertices, obtaining a cut (not
necessarily a bisection) with each $e(V_i) \leq \big( \frac{1}{4} +
\varepsilon \big) m$.  This is easily converted into a perfect bisection
by suitably distributing the $\geq \frac{n}{2}$ isolated vertices.
\hfill $\Box$

\section{Large bisection}

\label{sec:large-bisection}

Our eventual goal is to prove Theorem \ref{thm:bisection-min}, which
completely answers the question of judicious bisections under precise
minimum-degree conditions.  Using Theorem \ref{thm:bs-conjecture} as a
starting point, we see that the difficult case is when both the total
number of edges and the maximum degree are linear in $n$.  The natural
route is therefore to pre-allocate the vertices of large degree, and
apply randomized arguments on the remainder of the graph.  In this
section, we develop tools which we apply to the graph after the removal of large degree
vertices.  Our first objective is to establish a gain above the
trivial bound for the size of a bisection.  (This is the precursor to
a judicious bisection.)  Specifically, we prove Theorems
\ref{thm:tight} and \ref{thm:almost-bisection}.

The basic idea is to execute a deterministic greedy algorithm on a
suitably ordered input.  Indeed, suppose the vertices have been
partitioned into pairs $P_1, \ldots, P_{n/2}$.  We then obtain a
bisection $V = V_1 \cup V_2$ by separating each pair.  Start by
arbitrarily splitting $P_1$, and iterate as follows.  When considering
$P_i = \{v, w\}$, there are two ways to split: either we place $v \in
V_1$ and $w \in V_2$, or we place $v \in V_2$ and $w \in V_1$.  Let
$x_{i-1}$ be the number of crossing edges in the bisection induced by
the previously generated partition of $P_1 \cup \ldots \cup P_{i-1}$.
If $x_i$ is the number after splitting $P_i$, then an averaging
argument shows that there is a way to split $P_i$ such that
\begin{equation}
  \label{ineq:greedy-split}
  x_i - x_{i-1} \geq \frac{1}{2}
  \left[
    e(P_1 \cup \cdots \cup P_i)
    -
    e(P_1 \cup \cdots \cup P_{i-1})
  \right]
\end{equation}
Our algorithm proceeds by always selecting the split which maximizes
$x_i - x_{i-1}$.  Observe that this immediately implies that there is
a bisection of size at least $\frac{m}{2}$.  Our objective is now to
improve upon this bound by taking advantage of certain situations in
which \eqref{ineq:greedy-split} is not tight.

For example, if $P_i = \{v, w\}$ was in fact an edge, then both ways
of splitting $P_i$ would capture that edge.  Inequality
\eqref{ineq:greedy-split} only assumes that at least one of the splits
captures the edge, so in this case we could actually improve the bound
to
\begin{equation}
  \label{ineq:greedy-split-gain}
  x_i - x_{i-1} \geq \frac{1}{2}
  \left[
    e(P_1 \cup \cdots \cup P_i)
    -
    e(P_1 \cup \cdots \cup P_{i-1})
  \right]
  +
  \frac{1}{2}
\end{equation}

Alternatively, even if $P_i = \{v, w\}$ were not an edge, there is
another way to obtain the same bound as in
\eqref{ineq:greedy-split-gain}.  Indeed, whenever $e(P_1 \cup
\cdots \cup P_i) - e(P_1 \cup \cdots \cup P_{i-1})$ is an odd number,
it is impossible for \eqref{ineq:greedy-split} to be sharp, because
the right hand side is not an integer.  In such a situation, we
automatically obtain \eqref{ineq:greedy-split-gain}.  Note that
this happens when the ``back-degrees'', i.e., degrees of $v$ and
$w$ in $P_1 \cup \cdots \cup P_{i-1}$ are of different parity.

These two situations are advantageous because if
\eqref{ineq:greedy-split-gain} holds a linear number of times, then we
obtain our desired gain over the trivial bound by an amount linear in
$n$.

\subsection{Tight components and free vertices}

We aim to exploit the two possibilities outlined at the end of the
previous section.  In order to take full advantage of the first, we
begin by taking a matching of maximum size as the basis of our
pairing.  In this section, we demonstrate the connection between
maximum matchings and tight components.

Recall from the introduction (Definition \ref{def:tight}) that a tight
component is a connected component $T$ such that for every $v \in T$,
the subgraph induced by $T \setminus \{v\}$ contains a perfect
matching, and every perfect matching of $T \setminus \{v\}$ has the
property that no edge of the perfect matching has exactly one endpoint
adjacent to $v$.

\begin{definition}
  Let $\{e_1, \ldots, e_s\}$ be the edges of a maximum matching in a
  graph $G = (V, E)$, and let $W$ be the set of vertices not in the
  matching.  With respect to this fixed matching, say that a vertex
  $v$ in a matching edge $e_i$ is a \textbf{free neighbor} of a vertex
  $w \in W$ if $w$ is adjacent to $v$, but $w$ is not adjacent to the
  other endpoint of $e_i$.  Call a vertex $w \in W$ a \textbf{free
    vertex} if it has at least one free neighbor.
\end{definition}

The following lemma delineates the relationship between free vertices
and tight components.

\begin{lemma}
  \label{lem:freeandtight}
  Let $\{e_1, \ldots, e_s\}$ be the edges of a maximum matching in a
  graph $G = (V, E)$, and let $W$ be the set of vertices not in the
  matching.  Further assume that among all matchings of maximum size,
  we have chosen one which maximizes the number of free vertices in
  $W$.  Then, every tight component contains a distinct non-free
  vertex of $W$, and all non-free $W$-vertices are covered in this way
  (there is a bijective correspondence).
\end{lemma}

\noindent \textbf{Remark.}\, In particular, the number of non-free
vertices of $W$ is equal to the number of tight components of $G$.

\vspace{2mm}

\noindent \textbf{Proof.}\, The matching $\{e_1, \ldots, e_s\}$ must
be maximal within each connected component.  One basic property of a
tight component is that it contains an almost-perfect matching which
misses only one vertex.  Consequently, by maximality, in every tight
component $T$, $\{e_1, \ldots, e_s\}$ must miss exactly one vertex $w
\in W$.  Furthermore, the second property of a tight component is that
$w$ must have either 0 or 2 neighbors in each edge $e_i$ in $T$ (and
$w$ must have 0 neighbors in each edge $e_j$ not in $T$, since $T$ is
the connected component containing $w$).  Therefore, the unique vertex
$w$ is in fact a non-free $W$-vertex contained in $T$.

The remainder of the proof concentrates on the more substantial part
of the claim, which is that each non-free $W$-vertex is contained in
some tight component.  Consider such a vertex $w$, and let $T$ be a
maximal set of vertices which (i) contains $w$, (ii) induces a
connected graph which is a tight component, and (iii) does not cut any
$e_i$.  Since the set $\{w\}$ satisfies (i)--(iii), our optimum is taken
over a non-empty set, and so $T$ exists.

If $T$ is already disconnected from the rest of the graph, then we are
done.  So, consider a vertex $v \not \in T$ which has a neighbor $v'
\in T$.  If $v \in W$, then we can modify our matching by taking the
edge $vv'$, and changing the matching within $T$ by using property (ii) to
generate a new matching of $T \setminus \{v'\}$.  This will not affect
the matching outside of $T \cup \{v\}$, because property (iii) insulates the
adjustments within $T$ from the rest of the matching outside.  We
would then obtain a matching with one more edge, contradicting
maximality.  Therefore, all vertices $v \not \in T$ which have
neighbors in $T$ also satisfy $v \not \in W$.

Let us then consider a vertex $v_1 \not \in T$ with a neighbor $v' \in
T$.  We now know that $v_1$ must be covered by a matching edge; let
$v_2$ be the other endpoint of that edge.  By (iii), we also have $v_2
\not \in T$.  Note that $v_2$ cannot have a neighbor $w' \in W \setminus T$, or
else we could improve our matching by replacing $v_1 v_2$ with $w'
v_2$ and $v_1 v'$, and then using (ii) to take a perfect matching of
$T \setminus \{v'\}$.

Our next claim is that $v_2$ must be adjacent to $v'$ as well.
Indeed, assume for contradiction that this is not the case.  Then,
consider modifying our matching by replacing $v_1 v_2$ with the edge
$v_1 v'$ and changing the matching within $T$ by using (ii) to
generate a new matching of $T \setminus \{v'\}$.  As before, (iii)
ensures that the result is still a matching.  This time, the new
matching has the same size as the original one, but with more free
$W$-vertices (note that the vertex $v_2$ replaced the vertex $w$
in the set $W$). To see this, observe that $v_2$ is now unmatched and
free because it is adjacent to $v_1$ but not $v'$.
Previously, the only $W$-vertex inside $T$ was
our original $w$, which we assumed to be non-free in the first place.
Also, no other vertices outside of $T$ changed from being free to non-free,
because we already showed that no $W$-vertices outside of $T$ were
adjacent to $T \cup \{v_2\}$, and so any vertex that was free by
virtue of its adjacency with $v_1$ but not $v_2$ is still free because
it is not adjacent to $v'$ either.
This contradiction to maximality establishes that $v_2$ must be
adjacent to $v'$.

We now have $v_1$, $v_2$, and $v'$ all adjacent to each other, and no
vertices of $W \setminus T$ are adjacent to $T \cup \{v_1, v_2\}$.  Our argument
also shows that for any $v'' \in T$ which is adjacent to one of $v_1$
or $v_2$, it also must be adjacent to the other.  Our final objective
is to show that $T' = T \cup \{v_1, v_2\}$ also satisfies (i)--(iii),
which would contradict the maximality of $T$.  Properties (i) and
(iii) are immediate, so it remains to verify the conditions of a tight
component.  Since $T$ is tight and $v_1 v_2$ is an edge, $T' \setminus
\{u\}$ has a perfect matching for any $u \in T$.  The tightness of $T$
and the pairwise adjacency of $v_1$, $v_2$, and $v'$ also produce this
conclusion if $u \in \{v_1, v_2\}$.  It remains to show that for any
$u \in T'$ and any perfect matching of $T' \setminus \{u\}$, $u$ has
either 0 or 2 vertices in every matching edge.  But if this were not
the case, then we could replace the matching within $T'$ with the
violating matching of $T' \setminus \{u\}$.  The two matchings would
have the same size, but $u$ would become a free vertex.  No other
vertex of $W$ is adjacent to $T'$ by our observation above, so the
number of free vertices would increase, contradicting the maximality
of our initial matching.  Therefore, $T'$ induces a tight component,
contradicting the maximality of $T$.  We conclude that $T$ must have
been disconnected from the rest of the graph, as required.  \hfill
$\Box$

\subsection{Proof of Theorem \ref{thm:tight}}
\label{sec:thm-tight-proof}

In this section, we show that every $n$-vertex graph with $m$ edges,
maximum degree $\Delta$, and exactly $\tau$ tight components has a
bisection of size at least
\begin{displaymath}
  \frac{m}{2} + \frac{n - \max\{\tau, \Delta - 1\}}{4} \,.
\end{displaymath}
Our proof is somewhat similar to that used by
Erd\H{o}s, Gy\'arfas, and Kohayakawa \cite{ErGyKo97} to estimate Max Cut.

Start by taking a maximum matching $\{e_1, \ldots, e_s\}$ as in Lemma
\ref{lem:freeandtight}, which further maximizes the number of free
vertices in $W = V \setminus \{e_1, \ldots, e_s\} = \{w_1, \ldots,
w_r\}$.  Note that $n = r + 2s$.  By maximality, $W$ is an independent
set, and there cannot be any $e_i = \{v, v'\}$ such that $vw$ and
$v'w'$ are both edges, with $w, w'$ distinct vertices in $W$.  The
second observation shows that if two vertices $w, w' \in W$ each have
a free neighbor in an $e_i$, then their free neighbor is the same
vertex.

Let us now partition the vertices of $V$ into disjoint pairs.  Later,
we will distribute the vertices by separating each pair; this will
produce a bisection.  First, form a pair from each matching edge by
taking its endpoints.  In a slight abuse of notation, we will denote
these pairs by $e_1, \ldots, e_s$.  Then, for the vertices in $W$,
greedily form disjoint pairs by taking vertices $w_i$ and $w_j$ whose
sets of free neighbors are distinct.  (Note that if $w_i$ has no free
neighbors and $w_j$ has several, we may pair them.)  Let $p_1, \ldots,
p_{r'}$ be the pairs constructed from $W$ in this way.  Note that at
the end, we must have $r' \leq r/2$, and all of the $n-2(s+r')$
remaining vertices must have exactly the same set $S$ of free
neighbors.  Arbitrarily collect the remaining unpaired vertices into
pairs $q_1, \ldots, q_t$, where $t = \frac{n}{2} - (s+r')$.

If $S = \emptyset$, then all remaining vertices are non-free, so Lemma
\ref{lem:freeandtight} implies that $n-2(s+r') \leq \tau$, where
$\tau$ is the number of tight components.  Otherwise, there is some
vertex $v \in S$ which is adjacent to all remaining vertices.  Since
the degree of $v$ is at most $\Delta$ and one of its incident edges is
the matching edge which covers it, we find that $n-2(s+r') \leq \Delta
- 1$.  Therefore, $n-2(s+r') \leq \max\{ \tau, \Delta-1 \}$, and so it
suffices to show that there exists a bisection of size at least
$\frac{m}{2}+\frac{s+r'}{2}$ because
\begin{displaymath}
  \frac{m}{2} + \frac{s+r'}{2}
  \geq
  \frac{m}{2} + \frac{n - \max\{\tau, \Delta - 1\}}{4}
  \,.
\end{displaymath}

For this, we order the pairs.  Start by taking $e_1, e_2, \ldots, e_s$
in that order.  We iteratively insert the pairs $p_i$ into the
sequence as follows.  By construction, for each $p_i$, there exists an
edge $e_j$ such that only one of the vertices in $p_i$ has a free
neighbor in $e_j$.  Let $e_j$ be the lowest-indexed such edge and
place $p_{i}$ between $e_{j}$ and $e_{j+1}$ (arbitrarily order the
pairs $p_{i}$ which are placed between the same consecutive edges).
Finally, append the pairs $q_1, \ldots, q_t$ at the end of the
sequence, and call the resulting sequence $P_1, P_2, \ldots, P_{n/2}$.
(Here, each $P_i$ is either an $e_j$, a $p_j$, or a $q_j$.)

We split each pair across the bisection by using the $P_i$-sequence as
the input to the greedy algorithm presented at the beginning of
Section \ref{sec:large-bisection}.  By our discussion there, it
suffices to show that \eqref{ineq:greedy-split-gain} holds at least $s
+ r'$ times.  We already have $s$ pairs $P_i$ of the form $e_j$
(edges), which we know from the beginning to produce
\eqref{ineq:greedy-split-gain}.

We also have $r'$ pairs $P_i$ of the form $p_j = \{w, w'\}$.  There
are no edges from $w$ or $w'$ back to any vertices of $W$ in $P_1 \cup
\cdots \cup P_{i-1}$, because $W$ is an independent set.  As for edges
$e_j$ which appear among $P_1 \cup \cdots \cup P_{i-1}$, whenever $w$
(or $w'$) has exactly 1 edge to an endpoint of $e_j$, then that
endpoint is a free neighbor of $w$ (respectively. $w'$).  By construction, $w$
and $w'$ have exactly the same set of free neighbors in $P_1 \cup
\cdots \cup P_{i-1}$, except for one endpoint of an edge, at which
they differ.  Hence the total number of edges from $w$ and $w'$ back
to $P_1 \cup \cdots \cup P_{i-1}$ is odd, which we also showed at the
beginning to give \eqref{ineq:greedy-split-gain}.

Therefore, we indeed achieve \eqref{ineq:greedy-split-gain} at least
$s + r'$ times, and so the final bisection has size at least
$\frac{m}{2} + \frac{s+r'}{2}$, as required.  \hfill $\Box$

\subsection{Bisection corollaries}

We now use Theorem \ref{thm:tight} to prove Corollaries
\ref{cor:bisection+n/6-maxdegree} and \ref{cor:bisectionconnected}.
We begin by showing that every $n$-vertex graph with $m$ edges,
maximum degree at most $\frac{n}{3} + 1$, and no isolated vertices,
has a bisection of size at least $\frac{m}{2} + \frac{n}{6}$.

\vspace{2mm}

\noindent \textbf{Proof of Corollary
  \ref{cor:bisection+n/6-maxdegree}.}\, Tight components are always of
odd order, because the deletion of any vertex leaves a graph which has
a perfect matching.  Since we assume there are no isolated vertices,
this forces the number of tight components to be at most $n/3$, since
the smallest possible tight component is then of order 3.  Applying
Theorem \ref{thm:tight} with $\Delta \leq \frac{n}{3} + 1$ and $\tau
\leq \frac{n}{3}$, we find that $G$ has a bisection of size at least
\begin{displaymath}
  \frac{m}{2} + \frac{n - (n/3)}{4}
  =
  \frac{m}{2} + \frac{n}{6} \,,
\end{displaymath}
as desired.  \hfill $\Box$

\vspace{2mm}

Next, we show that every $n$-vertex connected graph with $m$ edges and
maximum degree $\Delta$ has a bisection of size at least $\frac{m}{2}
+ \frac{n + 1 - \Delta}{4}$.

\vspace{2mm}

\noindent \textbf{Proof of Corollary
\ref{cor:bisectionconnected}.}\, If $n=2$, then the graph consists
of a single edge, and the conclusion can easily be checked to be
true. Thus we may assume $n \ge 3$, from which we obtain $\Delta \ge
2$, since the graph is connected. Moreover, since the graph is
connected, the number of tight components is either 0 or 1. So
$\max\{\tau, \Delta - 1\} = \Delta - 1$. Therefore, Theorem
\ref{thm:tight} produces a bisection of size at least
\begin{displaymath}
  \frac{m}{2} + \frac{n - (\Delta - 1)}{4} \,,
\end{displaymath}
as desired.  \hfill $\Box$

\subsection{Approximate bisection}

Our general approach from Section \ref{sec:thm-tight-proof} turns out
to be quite powerful, as a slight twist of the approach produces
Theorem \ref{thm:almost-bisection}, which concerns $\alpha$-bisections
(bipartitions with both sides containing at least $(\frac{1}{2} -
\alpha)n$ vertices).  In this section, we prove that for any $0 \leq
\alpha \leq \frac{1}{6}$, every graph which has no isolated vertices
contains an $\alpha$-bisection of size at least $\frac{m}{2} + \alpha
n$.

\vspace{2mm}

\noindent \textbf{Proof of Theorem \ref{thm:almost-bisection}.}\, As
in Section \ref{sec:thm-tight-proof}, take an optimal matching $\{e_1,
\ldots, e_s\}$, greedily form the pairs $p_1, \ldots, p_{r'}$, and
insert them into the sequence of $e_i$'s as before.  However, this
time we will not necessarily pair up the remaining vertices.

\vspace{2mm}

\noindent \textbf{Case 1: There are at most $\boldsymbol{\frac{n}{3}}$
  unpaired vertices.}\, The number of unpaired vertices is exactly $n
- 2(s+r')$, so
\begin{align*}
  n - 2(s + r') &\leq \frac{n}{3} \\
  s + r' &\geq \frac{n}{3} \,.
\end{align*}
In this case, we now arbitrarily pair up the remaining vertices and
append the new pairs to the end of the sequence of pairs.  Applying
the same procedure and analysis as before, we obtain a
\emph{perfect}\/ bisection of size at least
\begin{displaymath}
  \frac{m}{2} + \frac{s+r'}{2}
  \geq
  \frac{m}{2} + \frac{n}{6}
  \geq
  \frac{m}{2} + \alpha n  \,,
\end{displaymath}
as desired.

\vspace{2mm}

\noindent \textbf{Case 2: There are more than
  $\boldsymbol{\frac{n}{3}}$ unpaired vertices.}\, Recall that the
remaining unpaired vertices are either all non-free, or all have the
same nonempty set of free neighbors.  Lemma \ref{lem:freeandtight},
bounds the number of non-free vertices by the number of tight
components.  Since all tight components are odd, and in this case we
have no isolated vertices, the number of tight components is at most
$\frac{n}{3}$.  Therefore, we must have that all unpaired vertices
have the same nonempty set of free neighbors.

Thus far, we have interleaved the $p_i$ into the original sequence of
$e_i$'s, and some vertices, say $w_1, \ldots, w_k$ with $k = n -
2(s+r') > \frac{n}{3}$, remain.  Importantly, note that both of
$\frac{k}{2} \pm \alpha n$ are between 0 and $k$ inclusive.

Let $e_i$ be the first matching edge which contains a free neighbor of
the remaining vertices.  Apply the greedy splitting algorithm
(beginning of Section \ref{sec:large-bisection}) to the sequence of
pairs up through $e_i$.  Let $v$ be the endpoint of $e_i$ that is
adjacent to all of $R = \{w_1, \ldots, w_k\}$.  Place $\frac{k}{2} +
\alpha n$ vertices of $R$ on the side opposite $v$, and the other
$\frac{k}{2} - \alpha n$ vertices of $R$ on the same side as $v$.
Then apply the greedy splitting algorithm to the rest of the pairs in
the sequence.

Note that by the trivial bound \eqref{ineq:greedy-split},
every split always adds at least half of the new edges to the cut, and
the cut we obtain in the end has parts of size
$(\frac{1}{2} \pm \alpha) n$. To finish the proof, we will show
that in the step when we split the vertices of $R$, there is an extra
gain of $\alpha n$ edges over the trivial bound of \eqref{ineq:greedy-split}.
For this, observe that
vertices of $R$ are adjacent to neither or both of the endpoints
of each $e_j$, $j < i$, but adjacent to $v \in e_i$.  By maximality of
the original matching, $W$ is an independent set, so there will be no
edges from $R$ to the pairs $p_j$ preceding $e_i$.  Let $b$ be the
number of edges $e_j$, $j < i$, which have both endpoints adjacent to
some vertex of $R$.

Regardless of how we split the vertices of $R$, exactly
$b$ edges out of the $2b$ edges that are incident to $e_j$ (for $j<i$)
will cross the cut.
If a vertex $w \in R$ is placed on the same side as $v \in
e_i$, then it will add no edge to the cut, and if it
is placed on the opposite side of
$v$, then it will add one edge to the cut.  Therefore,
by placing $\frac{k}{2} + \alpha n$ vertices of $R$ on the opposite
side of $v$, the total number of crossing edges from $R$ back to the
pairs preceding it is exactly
\begin{displaymath}
  b
  +
  \left( \frac{k}{2} + \alpha n \right)
%  +
%  \left( \frac{k}{2} - \alpha n \right) (b)
%  =
%  kb + \frac{k}{2} + \alpha n
  =
%  \frac{(2b + 1) k}{2} + \alpha n \,.
    \left(b + \frac{k}{2} \right) + \alpha n
\end{displaymath}
Since there were exactly $2b + k$ edges from $R$ back to the
preceding pairs, this produces the desired gain of $\alpha n$.
\hfill $\Box$

\section{Judicious bisection}
\label{sec:judicious-bisection}

Our previous results produced bipartitions which achieved certain cut
sizes.  In this section, we use the following martingale concentration
result (essentially the Hoeffding-Azuma inequality) to prove the
stronger Theorem \ref{thm:tightrandom}, which produces bisections that
induce few edges in each part.

\begin{theorem}
  \label{thm:azuma}
  (Corollary 2.27 in \cite{JLR}.)
  Given real numbers $\lambda, C_1, \ldots, C_n > 0$, let $f : \{0,
  1\}^n \rightarrow \mathbb{R}$ be a function satisfying the following
  Lipschitz condition: whenever two vectors $z, z' \in \{0, 1\}^n$
  differ only in the $i$-th coordinate (for any $i$), we always have
  $|f(z) - f(z')| \leq C_i$.  Suppose $X_1, X_2, \ldots, X_n$ are
  independent random variables, each taking values in $\{0, 1\}$.
  Then, the random variable $Y = f(X_1, \ldots, X_n)$ satisfies
  \begin{displaymath}
    \pr{ | Y - \E{Y} | \geq \lambda }
    \leq
    2 \exp\left\{
      - \frac{\lambda^2}{2 \sum C_i^2}
    \right\} \,.
  \end{displaymath}
\end{theorem}

\subsection{Core result}

We use this to control the performance of a randomized partitioning
algorithm, which will produce Theorem \ref{thm:tightrandom} (restated
here for the reader's convenience):

\begin{theorem*}
  Given any real constants $C, \varepsilon > 0$, there exist $\gamma, n_0
  > 0$ for which the following holds.  Every graph $G = (V, E)$ with
  $n \geq n_0$ vertices, $m \leq Cn$ edges, maximum degree at most
  $\gamma n$, and $\tau$ tight components, admits a bisection $V = V_1
  \cup V_2$ where each $V_i$ induces at most
  \begin{displaymath}
    \frac{m}{4} - \frac{n - \tau}{8} + \varepsilon n
  \end{displaymath}
  edges.
\end{theorem*}

\noindent \textbf{Proof.}\, Let us assume that $\varepsilon \leq
\frac{1}{2}$ and $C \geq 1$.  (Otherwise, we may reduce $\varepsilon$ to
$\frac{1}{2}$ or increase $C$ to 1, follow the proof below, and
achieve a stronger result.)  Let $\gamma = \frac{\varepsilon^4}{1024
  C^3}$.  We will implicitly assume large $n$ in the remainder of our
argument.  As in the proof of Theorem \ref{thm:tight}, start by taking
a maximum matching $\{e_1, \ldots, e_s\}$, breaking ties by counting
the number of free vertices in $W = V \setminus \{e_1, \ldots, e_s\} =
\{w_1, \ldots, w_r\}$.  Partition $W$ into $W_1 \cup \cdots \cup W_s
\cup T$ by placing all non-free vertices and vertices with degree
$\geq \frac{C}{\varepsilon}$ into $T$, and for each remaining $w \in W$,
place it in the lowest-indexed $W_i$ for which $w$ has a free-neighbor
in $e_i$.  For each $i$, let $U_i = e_i \cup W_i$.  Since there are
exactly $\tau$ non-free vertices by Lemma \ref{lem:freeandtight}, and
the sum of all degrees is at most $2Cn$, the combined size of $T$ is
at most
\begin{equation}
  \label{ineq:T-plus-big-degrees}
  |T| \leq \tau + 2 \varepsilon n \,.
\end{equation}

By our matching's maximality, no two vertices $w, w' \in W_i$ can be
adjacent to different endpoints of $e_i$.  So, all $w \in W_i$ are
only adjacent to the same vertex $v_i \in e_i$, and we conclude that
each $U_i$ induces a star with apex $v_i$.  We now randomly construct
a bipartition $V = V_1 \cup V_2$ (not yet a bisection) by
independently placing each $v_i$ on a uniformly random side, and then
placing the rest of $U_i \setminus \{v_i\}$ on the other side.
Independently place each vertex of $T$ on a uniformly random side.
Define the random variables $Y_1$ and $Y_2$ to be the numbers of edges
induced by each $V_i$.

The outline of the remainder of our proof is as follows.

\begin{description}
\item[Step 1.] The expectations $\E{Y_1}$ and $\E{Y_2}$ are both at
  most $\frac{m}{4} - \frac{n - \tau}{8} + \frac{\varepsilon n}{4}$.

\item[Step 2.] By applying the Hoeffding-Azuma concentration
  inequality (Theorem \ref{thm:azuma}) with suitable Lipschitz
  conditions, we find that with probability strictly greater than
  $\frac{2}{3}$, we have $Y_1 \leq \frac{m}{4} - \frac{n - \tau}{8} +
  \frac{\varepsilon n}{2}$.  By symmetry, the same probabilistic bound
  holds for $Y_2$.

\item[Step 3.] Another application of Hoeffding-Azuma shows that with
  probability greater than $\frac{2}{3}$, the sizes of $V_1$ and $V_2$
  are both within $\frac{\varepsilon n}{16 C}$ of $\frac{n}{2}$.

\item[Step 4.] There must therefore be an outcome in which all above
  properties hold simultaneously.  We now equalize the sizes of $V_1$
  and $V_2$ by moving at most $\frac{\varepsilon n}{16 C}$ vertices from
  the larger side to the smaller, and show that this only changes each
  $e(V_i)$ by at most $\frac{\varepsilon n}{2}$.
\end{description}

\noindent \emph{Step 1.}\, We need to bound $\E{Y_1}$ and establish
the Lipschitz condition.  By symmetry, note that $\E{Y_1} = \E{Y_2}$,
and so
\begin{equation}
  \label{eq:induced-from-crossing}
  \E{Y_1} = \frac{1}{2} \left(
    m - \E{e(V_1, V_2)}
  \right) \,.
\end{equation}
We calculate $e(V_1, V_2)$ using the accounting procedure introduced
at the beginning of Section \ref{sec:large-bisection}.  That is, we
let $X_i$ be the number of crossing edges induced by the splits of
$U_1 \cup \cdots \cup U_i$, and estimate each $X_i - X_{i-1}$.  This
is the number of new crossing edges that are determined after
splitting $U_i$.  To estimate $\E{X_i - X_{i-1}}$, let $k$ be the
number of edges from $U_i$ to previous $U_j$.

By symmetry, the expected number of crossing edges from $U_i$ to
earlier $U_j$ is exactly $k/2$.  The gain comes from the fact that all
of the $|U_i| - 1$ edges in the star induced by $U_i$ always become
crossing edges by construction.  Putting all of these calculations
together, we find that
\begin{align*}
  \E{X_i - X_{i-1}}
  &=
  \frac{k}{2} + \frac{|U_i| - 1}{2} + \frac{|U_i| - 1}{2} \\
  &=
  \frac{1}{2} \left[
    e(U_1 \cup \cdots \cup U_i)
    -
    e(U_1 \cup \cdots \cup U_{i-1})
  \right]
  +
  \frac{|U_i| - 1}{2} \,.
\end{align*}
Finally, since each vertex of $T$ is independently assigned at the
end, on expectation exactly half of the edges incident to $T$
become crossing edges.  Therefore,
\begin{displaymath}
  \E{e(V_1, V_2)}
  = \frac{m}{2} + \sum_{i=1}^s \frac{|U_i| - 1}{2}
  = \frac{m}{2} + \frac{n - |T| - s}{2}
\end{displaymath}
We can bound $s \leq \frac{n - |T|}{2}$ because each $U_i$ has at
least 2 vertices (from $e_i$).  Hence, using
\eqref{ineq:T-plus-big-degrees} to control $|T|$, we obtain
\begin{displaymath}
  \E{e(V_1, V_2)}
  \geq \frac{m}{2} + \frac{n - |T| - \frac{n-|T|}{2}}{2}
  \geq \frac{m}{2} + \frac{n - (\tau + 2 \varepsilon n)}{4}
  \,,
\end{displaymath}
which together with \eqref{eq:induced-from-crossing} yields
\begin{displaymath}
  \E{Y_1}
  \leq
  \frac{1}{2} \left(
    m - \left(
      \frac{m}{2} + \frac{n - \tau - 2\varepsilon n}{4}
    \right)
  \right)
  =
  \frac{m}{4} - \frac{n - \tau}{8} + \frac{\varepsilon n}{4}
  \,,
\end{displaymath}
as desired.

\vspace{2mm}

\noindent \emph{Step 2.}\, For each $1 \leq i \leq s$, let $C_i$ be
the sum of the degrees of all vertices in $U_i$.  Clearly, flipping
the assignment of $v_i$ cannot affect $Y_1$ by more than $C_i$.  Also,
flipping the assignment of any $w \in T$ cannot change $Y_1$ by more
than the degree $d(w)$ of $w$.  Therefore, if we define
\begin{displaymath}
  L = \sum_{i=1}^s \left(
    \sum_{u \in U_i} d(u)
  \right)^2
  +
  \sum_{w \in T} d(w)^2 \,,
\end{displaymath}
the Hoeffding-Azuma inequality (Theorem \ref{thm:azuma}) gives
\begin{equation}
  \label{ineq:Y1-with-L}
  \pr{Y_1 \geq \E{Y_1} + \frac{\varepsilon n}{4}}
  \leq
  2 \exp\left\{
    - \frac{\varepsilon^2 n^2}{32 L}
  \right\} \,.
\end{equation}
Let us now control $L$.  We observed in Step 1 that each $U_i$ induced
a star with apex $v_i \in e_i$. Let $u_i$ be the other endpoint of
$e_i$.  This gives $|U_i| \leq d(v_i) + 1$, and by construction, every
other vertex in $U_i$ has degree less than $\frac{C}{\varepsilon}$.  So,
\begin{displaymath}
  \sum_{u \in U_i} d(u)
  \leq
  d(v_i) + d(u_i) + (d(v_i) - 1) \frac{C}{\varepsilon}
  \leq
  (d(u_i) + d(v_i)) \frac{2C}{\varepsilon} \,,
\end{displaymath}
and hence
% (a+b)^2 vs 2a^2 + 2b^2
% 2ab     vs a^2 + b^2
% sqrt(ab) vs sqrt( (a^2 + b^2)/2 )
% GM      <= RMS
\begin{align*}
  L &\leq
  \frac{4C^2}{\varepsilon^2} \sum_{i=1}^s (d(u_i) + d(v_i))^2 + \sum_{w \in T} d(w)^2 \\
  &\leq
  \frac{8C^2}{\varepsilon^2} \sum_{i=1}^s (d(u_i)^2 + d(v_i)^2) + \sum_{w \in T} d(w)^2 \\
  &\leq
  \frac{8C^2}{\varepsilon^2} \sum_{v \in V} d(v)^2 \le \frac{8C^2}{\varepsilon^2} (\gamma n) \sum_{v \in V} d(v)\\
  &\leq
  \frac{8C^2}{\varepsilon^2} (\gamma n) (2 Cn) \,,
\end{align*}
where we used that all $d(v) \leq \gamma n$ and that sum of all
degrees in $G$ is at most $2 Cn$.  However, we chose $\gamma =
\frac{\varepsilon^4}{1024 C^3}$, so we have $L \leq
\frac{\varepsilon^2 n^2}{64}$.  Substituting this into
\eqref{ineq:Y1-with-L}, we conclude that
\begin{displaymath}
  \pr{Y_1 \geq \E{Y_1} + \frac{\varepsilon n}{4}}
  \leq
  2 e^{-2}
  <
  \frac{1}{3} \,,
\end{displaymath}
as desired.

\vspace{2mm}

\noindent \emph{Step 3.}\, By symmetry, the expected size of $V_1$ is
exactly $\frac{n}{2}$.  Switching the choice of a single $v_i$ can
only change $|V_1|$ by at most $|U_i| - 2 \le d(v_i) - 1$, and switching the
choice of a single $w \in T$ can only change $|V_1|$ by at most 1.
Therefore, if we define
\begin{displaymath}
  L' = \sum_{i=1}^s ( d(v_i) - 1 )^2 + |T| \,,
\end{displaymath}
the Hoeffding-Azuma inequality (Theorem \ref{thm:azuma}) gives
\begin{equation}
  \label{ineq:V1-with-L}
  \pr{ \left| |V_1| - \frac{n}{2} \right| \geq \frac{\varepsilon n}{16 C} }
  \leq
  2 \exp\left\{
    - \frac{\varepsilon^2 n^2}{512 C^2 L'}
  \right\} \,.
\end{equation}
As in Step 2, we may control $L'$ by
\begin{align*}
  L'
  &\leq |T| + \sum_{i=1}^s d(v_i)^2 
  \leq n + (\gamma n) \sum_{i=1}^s d(v_i) 
  \leq n + (\gamma n) (2Cn) 
  = (1 + o(1)) \frac{\varepsilon^4 n^2}{512 C^2}
\end{align*}
Substituting this into \eqref{ineq:V1-with-L}, we find that
\begin{displaymath}
  \pr{ \left| |V_1| - \frac{n}{2} \right| \geq \frac{\varepsilon n}{16 C} }
  \leq
  2 e^{ -(1 + o(1)) / \varepsilon^2 }
  < \frac{1}{3} \,.
\end{displaymath}
Since $|V_1| + |V_2| = n$, whenever $|V_1|$ is within $\frac{\varepsilon
  n}{16C}$ of $\frac{n}{2}$, $|V_2|$ is as well.  Hence with
probability strictly greater than $\frac{2}{3}$, both $|V_1|$ 
and $|V_2|$ are within $\frac{\varepsilon n}{16 C}$ of $\frac{n}{2}$.

\vspace{2mm}

\noindent \emph{Step 4.}\, Our probabilistic arguments show that there
is a partition $V = V_1 \cup V_2$ with the properties that each
$|V_i|$ is within $\frac{\varepsilon n}{16 C}$ of $\frac{n}{2}$, and each
$e(V_i)$ is at most $\frac{m}{4} - \frac{n-\tau}{8} + \frac{\varepsilon
  n}{2}$.  To equalize the sizes of $|V_1|$ and $|V_2|$, we now move
the lowest-degree vertices from the larger side to the smaller.

Since the sum of all degrees is at most $2Cn$, the number of vertices
whose degree exceeds $8C$ is at most $\frac{n}{4}$.  In particular, we
may equalize $|V_1|$ and $|V_2|$ by moving vertices of degree below
$8C$.  This affects each $e(V_i)$ by at most $8C \cdot \frac{\varepsilon
  n}{16 C} = \frac{\varepsilon n}{2}$, thereby producing the desired
result.  \hfill $\Box$

\subsection{Judicious bisection with initial pre-partition}

Theorem \ref{thm:tightrandom} produces a judicious bisection which in
fact surpasses the expected performance of the elementary random
bisection algorithm, when degrees are bounded by $o(n)$ and there are
only $O(n)$ edges.  In Section \ref{sec:bs-conjecture}, we showed that
when there are $\omega(n)$ edges, then the elementary random bisection
algorithm already achieves the bound of Theorem \ref{thm:bisection-min}.
 So, the main obstacle is the maximum degree condition.

To work around this, we will often pre-allocate vertices of high
degree into $V_1$ and $V_2$, and then apply a randomized algorithm to
distribute the remaining vertices.  The same proof as we used in the
previous section will also produce the following corollary.

\begin{corollary}
  \label{cor:tightrandominit}
  Given any real constants $C, \varepsilon > 0$, there exist $\gamma, n_0
  > 0$ for which the following holds.  Let $G = (V, E)$ be a given
  graph with $n \geq n_0$ vertices and at most $Cn$ edges, and let $A
  \subset V$ be a set of $\leq \gamma n$ vertices which have already
  been partitioned into $A_1 \cup A_2$.  Let $\overline{A} = V
  \setminus A$, and suppose that every vertex in $\overline{A}$ has
  degree at most $\gamma n$ (with respect to the full $G$).  Let
  $\tau$ be the number of tight components in $G[\overline{A}]$.
  Then, there is a bisection $V = V_1 \cup V_2$ with $A_1 \subset V_1$
  and $A_2 \subset V_2$, such that both
  \begin{align*}
    e(V_1) &\leq e(A_1) + \frac{e(A_1, \overline{A})}{2}
    + \frac{e(\overline{A})}{4} - \frac{n - \tau}{8} + \varepsilon n \\
    e(V_2) &\leq e(A_2) + \frac{e(A_2, \overline{A})}{2}
    + \frac{e(\overline{A})}{4} - \frac{n - \tau}{8} + \varepsilon n \,.
  \end{align*}
\end{corollary}

\noindent \textbf{Proof.}\, Again, assume that $\varepsilon \leq
\frac{1}{2}$ and $C \geq 1$.  Let $\gamma = \frac{\varepsilon^4}{1024
  C^3}$.  We will implicitly assume large $n$ in the remainder of our
argument.  As in the proof of Theorem \ref{thm:tightrandom}, start by
taking a maximum matching $\{e_1, \ldots, e_s\}$ in $G[\overline{A}]$,
breaking ties by counting the number of free vertices in $W =
\overline{A} \setminus \{e_1, \ldots, e_s\} = \{w_1, \ldots, w_r\}$.
Partition $W$ into $W_1 \cup \cdots \cup W_s \cup T$ by placing all
non-free vertices and vertices with degree $\geq \frac{C}{\varepsilon}$
(with respect to the full $G$) into $T$, and for each remaining $w \in
W$, place it in the lowest-indexed $W_i$ for which $w$ has a
free-neighbor in $e_i$.  For each $i$, let $U_i = e_i \cup W_i$.  Note
that Lemma \ref{lem:freeandtight} bounds the number of non-free
vertices by the number of tight components $\tau$, and since the sum
of all degrees is at most $2Cn$, the combined size of $T$ is at most
\begin{displaymath}
  |T| \leq \tau + 2 \varepsilon n \,.
\end{displaymath}
As before, each $U_i$ induces a star with apex $v_i$.  Randomly
construct a bipartition by independently placing each $v_i$ on a
uniformly random side, and then placing the rest of $U_i \setminus
\{v_i\}$ on the other side.  Independently place each vertex of $T$ on
a uniformly random side.  Define the random variables $Y_1$ and $Y_2$
to be the numbers of edges induced by each $V_i$.  A similar argument
to that of Step 1 before now establishes that
\begin{align*}
  \E{Y_1} &= e(A_1) + \frac{e(A_1, \overline{A})}{2}
  + \frac{e(\overline{A})}{4} - \frac{n - \tau}{8} + \frac{\varepsilon n}{4} \\
  \E{Y_2} &= e(A_2) + \frac{e(A_2, \overline{A})}{2}
  + \frac{e(\overline{A})}{4} - \frac{n - \tau}{8} + \frac{\varepsilon n}{4} \,.
\end{align*}
The same considerations as in Step 2 before then establish that each
\begin{align*}
  \pr{Y_i \geq \E{Y_i} + \frac{\varepsilon n}{4}}
  &\leq
  2 e^{-2}
  <
  \frac{1}{3} \\
  \pr{ \left| |V_1 \setminus A_1| - \frac{|\overline{A}|}{2} \right| \geq \frac{\varepsilon n}{16\sqrt{2} C} }
  &\leq
  2 e^{ -1 / (2\varepsilon^2) }
  \le 2 e^{-2}
  < \frac{1}{3} \,,
\end{align*}
so there is an outcome where each $|V_i \setminus A_i|$ deviates from
its mean by at most $\frac{\varepsilon n}{16\sqrt{2} C}$, and each
\begin{displaymath}
  e(V_i) \leq e(A_i) + \frac{e(A_i, \overline{A})}{2}
  + \frac{e(\overline{A})}{4} - \frac{n - \tau}{8} + \frac{\varepsilon n}{2} \,.
\end{displaymath}
Since $|A| \leq \gamma n$, we can now equalize the sizes of $V_1$ and
$V_2$ by moving fewer than $\frac{\varepsilon n}{16\sqrt{2}C} + |A| \le \frac{\varepsilon n}{16C}$ of the
lowest-degree vertices from the larger side to the smaller side.  The
same argument as in Step 4 before then completes the proof.
\hfill $\Box$

\section{Minimum degree two}
\label{sec:bisection-min-2}

Corollary \ref{cor:tightrandominit}, our final result from the
previous section, now enables us to prove Theorem
\ref{thm:bisection-min}, which achieves the asymptotically optimal
bisection in graphs with prescribed minimum degree $\delta$.  Since
the arguments become increasingly complicated as $\delta$ gets larger,
we start with the $\delta = 2$ case. The proof of this
simple case already contains (in a more simple form) most of the
main ideas for the proof of the general case, and thus we hope that
its proof serves as an illustrative starting point. In this section,
we prove the following result.

\begin{theorem*}
  For every $\varepsilon > 0$, there exists $n_0$ such that every graph
  $G = (V, E)$ with $n > n_0$ vertices and minimum degree at least 2
  has a bisection $V = V_1 \cup V_2$ in which each $e(V_i) \leq \big(
  \frac{1}{3} + \varepsilon \big) e(G)$.
\end{theorem*}

\noindent \textbf{Proof.}\, Let $m$ be the number of edges in $G$.  By
Theorem \ref{thm:bs-conjecture} (with $\varepsilon = \frac{1}{12}$),
and by the final argument in its proof,
we may assume that $\frac{n}{4} \leq m \leq 144n$.
This observation
translates into $m = \Theta(n)$, where the notation has masked
absolute constants $\frac{1}{4}$ and $144$.

If all degrees are below $n^{3/4}$, then we are already done by
Theorem \ref{thm:tightrandom}.  Otherwise, let $A$ be the set of
vertices with degree at least $n^{3/4}$ and let $\overline{A} = V
\setminus A$.  Since $m = \Theta(n)$, we have $|A| = O(n^{1/4})$, and
hence $e(A) = o(m)$.  Let $m_1 = e(A, \overline{A})$ and $m_2 =
e(\overline{A})$, and note that $m = m_1 + m_2 + o(m)$.  Since $e(A) =
o(m)$, we will essentially ignore the edges inside $A$.  Let $\Delta$
be the maximum $\overline{A}$-degree in $A$, i.e., the maximum number
of edges to $\overline{A}$ that any vertex of $A$ has.

Start the bipartition $V = V_1 \cup V_2$ by splitting $A = A_1 \cup
A_2$ in such a way that $e(A_1, \overline{A}) \geq e(A_2,
\overline{A})$, but $e(A_1, \overline{A})$ and $e(A_2, \overline{A})$
are as close as possible. By the optimality of the splitting,
we have $e(A_1, \overline{A}) - e(A_2, \overline{A}) \le \Delta$, since otherwise
we can improve this partition by moving one vertex from $A_1$ to $A_2$.
Note that if this achieves $\frac{m_1}{3}
\leq e(A_1, \overline{A}) \leq \frac{2m_1}{3}$, then $e(A_2,
\overline{A})$ also satisfies the same bound, and Corollary
\ref{cor:tightrandominit} applied with $\tau \le n$
then produces a bisection with the desired
\begin{displaymath}
  e(V_i) \leq e(A) + \frac{2 m_1 / 3}{2} + \frac{m_2}{4} + o(m)
  \leq
  \frac{m}{3} + o(m)
\end{displaymath}

Therefore, we may assume that the vertices of $A$ cannot be split as
evenly as above.  Yet if $\frac{3}{2} \Delta \le m_1 \le 3\Delta$,
then placing the vertex with $\Delta$ on one side, and the rest of $A$
on the other already gives $\frac{m_1}{3} \leq e(A_1, \overline{A})
\leq \frac{2m_1}{3}$.  Also, if $m_1 > 3\Delta$, then the optimal
split of $A$ will give
\begin{displaymath}
  \frac{m_1}{2}
  \leq
  d(A_1, \overline{A})
  \leq
  \frac{m_1}{2} + \frac{\Delta}{2} < \frac{2m_1}{3} \,.
\end{displaymath}
Consequently, it remains to study the situation when $\Delta \leq m_1
< \frac{3}{2} \Delta$.  The optimal split of $A$ then consists of
putting the $\Delta$-vertex alone in $A_1$, and placing all other
vertices of $A$ into $A_2$.

\vspace{2mm}

\noindent \textbf{Case 1: $\boldsymbol{m_2 \geq 6 \Delta - 4m_1}$.}\,
We now have $e(A_1, \overline{A}) = \Delta$ and $e(A_2, \overline{A})
= m_1 - \Delta < \frac{\Delta}{2}$.  Corollary
\ref{cor:tightrandominit} applied with $\tau \le n$
then produces a bisection with the desired
\begin{displaymath}
  e(V_i) \leq e(A) + \frac{\Delta}{2} + \frac{m_2}{4} + o(m)
  \leq
  \frac{4 m_1 + m_2}{12} + \frac{m_2}{4} + o(m)
  =
  \frac{m}{3} + o(m) \,.
\end{displaymath}

\vspace{2mm}

\noindent \textbf{Case 2: $\boldsymbol{m_2 < 6 \Delta - 4m_1}$.}\,
To settle this case, we will first bound the number of tight
components in terms of $n$ and $m_2$. Let $n' = |\overline{A}|$, and
let $\tau_1$ be the number of isolated vertices in $G[\overline{A}]$.
Since $A_1$ consists of a single vertex, there are
exactly $m_1 - \Delta$ edges from $A_2$ to $\overline{A}$.  The
graph $G$ had minimum degree at least 2, so every isolated vertex in
$G[\overline{A}]$ consumes at least one edge from $A_2$ (it can only
lose at most one edge to the lone vertex in $A_1$). Therefore, $\tau_1
\le m_1 - \Delta$.

Let $t$ be the maximum number of vertex-disjoint triangles in $G[\overline{A}]$.
Since $G[\overline{A}]$ has $\tau_1$ isolated vertices, we immediately have $m_2 \geq
3t + \frac{n' - \tau_1 - 3t}{2}$, or equivalently, $3t \leq 2m_2 - n' + \tau_1$.
Every tight component of $G[\overline{A}]$ either is an isolated vertex,
or contains a triangle, so we may bound the
number of tight components $\tau$ by
\[ \tau \le \tau_1 + t \le \tau_1 + \frac{2m_2 - n' + \tau_1}{3} = \frac{4\tau_1 + 2m_2 - n'}{3}\,. \]
By construction, all
degrees in $G[\overline{A}]$ are below $n^{3/4}$, so Theorem \ref{thm:tight}
produces a bisection $\overline{A} = B_1 \cup B_2$ of $G[\overline{A}]$ of size at least
\begin{displaymath}
  \frac{m_2}{2} + \frac{n' - \tau - n^{3/4}}{4}
  \ge
  o(m) + \frac{m_2}{2} + \frac{n'}{4} - \frac{4\tau_1 + 2m_2 - n'}{12} = o(m) + \frac{m_2 + n' - \tau_1}{3}\,,
\end{displaymath}
which by $\tau_1 \le m_1 - \Delta$ is at least
\begin{displaymath}
  e(B_1, B_2) \ge o(m) + \frac{m_2 + n' - m_1 + \Delta}{3}\,.
\end{displaymath}
Without loss of generality, suppose that $e(A_1, B_1) + e(B_1)
\leq e(A_1, B_2) + e(B_2)$.  (Otherwise, we may swap $B_1$ and $B_2$.)
Let $V_i = A_i \cup B_i$.  Then, we have
\begin{align*}
  e(V_1)
  =
  e(A_1) + e(A_1, B_1) + e(B_1)
  &\leq
  e(A_1) + \frac{e(A_1, B_1) + e(A_1, B_2) + e(B_1) + e(B_2)}{2} \\
  &=
  o(m) + \frac{\Delta + (m_2 - e(B_1, B_2))}{2} \\
  &\leq
  o(m) + \frac{\Delta + m_2}{2} - \frac{m_2 + n' - m_1 + \Delta}{6} \\
  &=
  o(m) + \frac{\Delta}{3} + \frac{m_2}{3} - \frac{n'}{6} + \frac{m_1}{6} \,.
\end{align*}
Using $n' \ge \Delta$, $\Delta \le m_1$, and $m = m_1 + m_2 + o(m)$, this is at most
\begin{align*}
  e(V_1) &\leq
  o(m) + \frac{\Delta}{6} + \frac{m_2}{3} + \frac{m_1}{6}
  \le
  o(m) + \frac{m_1}{3} + \frac{m_2}{3} = o(m) + \frac{m}{3} \,.
\end{align*}
On the other hand,
\begin{align*}
  e(V_2) = e(A_2) + e(A_2, B_2) + e(B_2)
  &\leq
  o(m) + e(A_2, \overline{A}) + \left( e(\overline{A}) - e(B_1, B_2) \right) \\
  &\leq
  o(m) + (m_1 - \Delta) + \left(
    m_2 - \frac{m_2 + n' - m_1 + \Delta}{3}
  \right) \\
  &=
  o(m) + \frac{4m_1 - 4\Delta + 2m_2 - n'}{3} \\
  &=
  o(m) + \frac{m_1 + m_2}{3} + \frac{3m_1 - 4\Delta + m_2 - n'}{3} \,.
\end{align*}
Using our
assumption that $m_2 < 6\Delta - 4m_1$ in this case, this is at most
\begin{align*}
  e(V_2)
  &\le o(m) + \frac{m}{3} + \frac{2\Delta - m_1 - n'}{3}\,.
\end{align*}
Since $\Delta \leq m_1$ and $\Delta \le n' + o(n)$, we conclude that $e(V_2)
\leq o(m) + \frac{m}{3}$, as desired.

Therefore, our partition $V = V_1 \cup V_2$ has both $e(V_i) \leq
\frac{m}{3} + o(m)$.  Also, $|V_1|$ and $|V_2|$ differ by at most $|A|
= O(n^{1/4})$.  Since $m \leq 144n$, at most $\frac{n}{2}$
vertices of $G$ can have degree over $576$. Hence, there are
at least $\frac{n}{2}$ vertices with degree at most $576$.
We equalize $|V_1|$ and $|V_2|$ by moving $O(n^{1/4})$ such vertices across.  This
affects each $e(V_i)$ by under $O(n^{1/4}) = o(m)$, thereby completing
the proof.  \hfill $\Box$

\section{Minimum degree: general case}
\label{sec:bisection-min}

In this section,
we establish the general case of our main theorem, which we restate
below for the reader's convenience.

\begin{theorem*}
  Given any $\varepsilon > 0$ and any positive even integer $\delta$,
  there exists $n_0$ such that every graph $G = (V, E)$ with $n \geq
  n_0$ vertices, $m$ edges, and minimum degree at least $\delta$, has
  a bisection in which for both $i=1,2$,
  \begin{displaymath}
    e(V_i) \leq \left( \frac{\delta + 2}{4 (\delta + 1) } + \varepsilon \right) m \,.
  \end{displaymath}
\end{theorem*}

\noindent \textbf{Proof.}\, By Theorem \ref{thm:bs-conjecture} and the
final argument in its proof, we may assume that $\frac{n}{4} \leq m
\leq \varepsilon^{-2} n$.  Throughout this proof, $\varepsilon$ will be
fixed, and it will be convenient to employ asymptotic notation which
implicitly assumes that $\varepsilon$ is a constant.  For example, our
initial observation translates into $m = \Theta(n)$.

We first describe the idea of the remainder of the proof, which is similar
to that of the case $\delta = 2$ given in the previous section.
Start by identifying the set of
vertices $A$ which have large degree (degree at least $n^{3/4}$), and
let $A_1 \cup A_2$ be a bipartition of these vertices so that
$e(A_1, V \setminus A)$ and $e(A_2, V \setminus A)$ are as close as possible.
If the gap between these two is small enough, then a naive random
split of the remaining vertices will give the result we were hoping for.
Otherwise, we would like to make use of Corollary~\ref{cor:tightrandominit}.
Most of the technical details of our proof lie in this step, or more
precisely, in estimating the number of tight components in $G[V \setminus A]$.
Recall that to do this in the $\delta=2$ case, we first deduced a structural
information about the degree of the vertices in $A$ (there we
knew that $A_1$ consisted of a single vertex), and then used
this information to effectively bound the number of tight components.
The general idea is the same here, but things get more complicated
because both the structural information, and the method of
bounding tight components under this information becomes a lot more difficult
(and computationally involved).

We now give the details of the proof.
If all degrees are below $n^{3/4}$, then we are already done by
Theorem \ref{thm:tightrandom}.  Otherwise, let $A$ be the set of
vertices with degree at least $n^{3/4}$ and let $\overline{A} = V
\setminus A$ (we call the vertices in $A$ as the {\em large degree}
vertices).  Since $m = \Theta(n)$, we have $|A| = O(n^{1/4})$, and
hence $e(A) = o(m)$.  Let
\begin{displaymath}
  n' = |\overline{A}| \,,
  \quad\quad
  m_1 = e(A, \overline{A}) \,,
  \quad \text{and} \quad
  m_2 = e(\overline{A}) \,,
\end{displaymath}
and note that
\begin{equation}
  \label{eq:n-m-close}
  n' = n - o(n)
  \quad \text{and} \quad
  m = m_1 + m_2 + o(m) \,.
\end{equation}
Start the bipartition $V = V_1 \cup V_2$ by splitting $A = A_1 \cup
A_2$ in such a way that $e(A_1, \overline{A}) \geq e(A_2,
\overline{A})$, but $e(A_1, \overline{A})$ and $e(A_2, \overline{A})$
are as close as possible.  Define the parameter
\begin{displaymath}
  \theta = e(A_1, \overline{A}) - e(A_2, \overline{A})  \,,
\end{displaymath}
which quantifies the ``gap'' in the optimal split.  By construction,
we now have
\begin{displaymath}
  e(A_2, \overline{A}) \leq e(A_1, \overline{A})
  =
  \frac{m_1 + \theta}{2} \,,
\end{displaymath}
so Corollary \ref{cor:tightrandominit} produces a bisection in which
both
\begin{displaymath}
  e(V_i)
  \leq
  \frac{m_1 + \theta}{4} + \frac{m_2}{4} - \frac{n - \tau}{8} + o(m)
  \,,
\end{displaymath}
where $\tau$ is the number of tight components in $G[\overline{A}]$.
Using \eqref{eq:n-m-close} to combine $m_1$ and $m_2$, we observe that
it therefore suffices to show that
\begin{displaymath}
  \frac{m}{4} + \frac{\theta}{4} - \frac{n - \tau}{8}
  \leq
  \frac{\delta + 2}{4(\delta + 1)} m \,,
\end{displaymath}
or equivalently,
\begin{equation}
  \label{ineq:NEED:main}
  \theta + \frac{\tau}{2} \leq \frac{n}{2} + \frac{m}{\delta + 1}
  \,.
\end{equation}
Since we always have $\tau \leq n$, we see that we are already done if
$\theta \leq \frac{m}{\delta + 1}$.  (It is worth noting that $\theta
= \frac{m}{\delta + 1}$ corresponds to the extremal example
$K_{\delta+1, n - \delta - 1}$, as there $m \approx (\delta+1)n$, there are
exactly $\delta+1$ large degree vertices, and the optimal split puts
$\frac{\delta}{2}+1$ of them in $A_1$ and the other $\frac{\delta}{2}$ in $A_2$.)
Therefore, we will assume for the remainder of the proof that the
``gap'' between $e(A_1, \overline{A})$ and $e(A_2, \overline{A})$ is
rather large:
\begin{equation}
  \label{ineq:gap-at-least}
  \theta
  >
  \frac{m}{\delta + 1} \,.
\end{equation}

In fact, this gap is so wide that it has serious implications
regarding the nature of the $A = A_1 \cup A_2$ split.  Recall that in
the $\delta = 2$ case (previous section), we deduced that $A_1$ only
contained a single vertex of degree $\Delta$, and $e(A_2,
\overline{A})$ was below $\frac{\Delta}{2}$.  The following result,
which we prove in Section \ref{sec:splitoflargevertices}, is its
analogue in the general case.  Remarkably, its proof only uses the
optimality of the $A = A_1 \cup A_2$ split.  In order to maintain
continuity in our exposition, we postpone the proof of this lemma to
Section \ref{sec:splitoflargevertices} below.

\begin{lemma}
  \label{lem:splitoflargevertices2}
  If $\theta > \frac{m}{\delta + 1}$, then in the optimal split
  $A = A_1 \cup A_2$:
  \begin{description}
  \item[(i)] All $\overline{A}$-degrees in $A_1$ are at least $\theta$.
  \item[(ii)] Only at most $\delta - 1$ vertices in $A_1 \cup A_2$
    have $\overline{A}$-degree at least $\theta$.
  \item[(iii)] The sum of the $\overline{A}$-degrees over all other
    vertices (in $A_2$) is at most $n' - \theta$.
  \end{description}
\end{lemma}

\vspace{2mm}

\noindent \textbf{Remark.}\, We now also have $\theta \leq n'$,
because the $\overline{A}$-degree of the vertex in $A_1$ is at least
$\theta$, but cannot exceed $|\overline{A}| = n'$. We refer to the
vertices of $\overline{A}$-degree at least $\theta$ as the
{\em huge $\overline{A}$-degree vertices} or simply {\em huge degree vertices}.

\vspace{2mm}

Let $\alpha$ be the number of vertices in $A_1 \cup A_2$ with
$\overline{A}$-degree at least $\theta$, and let $\rho$ be the sum of
the $\overline{A}$-degrees of all other vertices (which are all in
$A_2$).  Note that the sum of the $\overline{A}$-degrees of the
vertices in $A$ is then
\begin{equation}
  \label{eq:m1-alpha}
  m_1 \ge \alpha \theta + \rho \,.
\end{equation}

Returning to verify \eqref{ineq:NEED:main}, we will need to control
the number of tight components. In the simple case when $\delta=2$,
we were able to do this by counting the number of isolated vertices,
and the number of vertex-disjoint triangles. However, things
get more difficult to control for general $\delta$, and
the most efficient approach is to
bound $\tau$ in terms of the degree sequence of $G[\overline{A}]$, via
the following result, which essentially claims that the optimum is
achieved when all vertices are packed into tight components that are
complete graphs of odd order.

\begin{lemma}
  \label{lem:tau-by-degree-sequence}
  Let $H$ be an arbitrary graph.  For each integer $i$, let $d_i$ be
  the number of vertices in $H$ with degree equal to $i$.  Then the
  number of tight components $\tau$ in $H$ satisfies
  \begin{displaymath}
    \tau
    \leq
    \frac{d_0}{1} + \frac{d_2}{3} + \frac{d_4}{5} + \cdots
  \end{displaymath}
\end{lemma}

\noindent \textbf{Remark.}\, Odd-degree vertices cannot contribute to
tight components because for every vertex $v$ in a tight component
$T$, there is a perfect matching of $T \setminus \{v\}$, and $v$ has
either 0 or 2 neighbors in each edge of the matching.  This is why the
bound above does not involve $d_1$, $d_3$, $d_5$, etc.

\vspace{2mm}

The degree sequence does not appear in \eqref{ineq:NEED:main},
however, and there are two ways to proceed and relate the degree sequence with
other parameters.  The simpler approach crudely
controls the degree sequence in terms of $\alpha$ and $m_1$ (or more
precisely $\alpha$ and $\rho$), using
the minimum degree condition of $G$.  (For example, if $m_1$ were 0,
then the minimum degree condition would force $d_0 = d_1 = \cdots =
d_{\delta - 1} = 0$.)  The more intricate approach adds an additional
information by relating the degree sequence with $m_2$ as well.

We will actually take both approaches by the end of this proof, but it
turns out to be more convenient to use the crude approach first,
because it will already work for most of the (non-tight) cases.
Applying the minimum degree condition in $G$ to the sum of all degrees
in $\overline{A}$, we find that
\begin{displaymath}
  \delta n'
  \leq
  2 m_2 + m_1 \,,
\end{displaymath}
which means that
\begin{equation}
  \label{ineq:m-by-delta-m1}
  m
  \geq
  m_1 + m_2
  \geq
  \frac{\delta n' + m_1}{2} \,.
\end{equation}
%Also, observe that even if all $\alpha$ ``high-$\overline{A}$-degree''
%vertices of $A$ had full degree to $\overline{A}$, the cumulative
%deficiency of the degrees in $G[\overline{A}]$ (with respect to the
%threshold $\delta - \alpha$) is at most $\rho$.
We then make use of the fact that $A$ contains $\alpha$
huge degree vertices.
Note that even if all the $\alpha$ vertices in $A$ had
full degree to $\overline{A}$, the only way that a vertex in
$G[\overline{A}]$ could have degree below $\delta - \alpha$ would be if it
was incident to edges counted by $\rho$. Consequently, we
get the following constraint:
\begin{equation}
  \label{ineq:degrees-deficiency}
  (\delta - \alpha) d_0
  + (\delta - \alpha - 1) d_1
  + \cdots
  + d_{\delta - \alpha - 1}
  \leq \rho \,.
\end{equation}

Then, it turns out that the most economical way to control the bound
in Lemma \ref{lem:tau-by-degree-sequence} is to use $\rho$ to reduce
as many degrees to 0 as possible.  This produces the following result,
whose complete proof is in Section \ref{sec:tau-by-rho}.

\begin{lemma}
  \label{lem:tau-by-rho}
  The induced subgraph $G[\overline{A}]$ has at most $\frac{n' +
    \rho}{\delta - \alpha + 1}$ tight components.
\end{lemma}

By \eqref{ineq:m-by-delta-m1} and \eqref{eq:m1-alpha},
we get $m \ge \frac{\delta n' + \alpha \theta + \rho}{2}$.
Using this bound and Lemma \ref{lem:tau-by-rho} in
\eqref{ineq:NEED:main}, we see that it suffices to show
\begin{displaymath}
  \theta + \frac{n' + \rho}{2(\delta - \alpha + 1)}
  \leq
  \frac{n}{2}
  + \frac{\delta n' + \alpha \theta + \rho}{2(\delta + 1)}
  \,,
\end{displaymath}
or since $n' \le n$, that
\begin{equation}
  \label{ineq:NEED:alpha-big}
  \theta + \frac{n' + \rho}{2(\delta - \alpha + 1)}
  \leq
  \frac{2\delta + 1}{2(\delta+1)} n'
  + \frac{\alpha \theta + \rho}{2(\delta + 1)}  \,.
\end{equation}
It turns out that the above inequality is true for most values of $\alpha$ and
$\delta$ (for the remaining cases, we will use the slightly different approach that
we briefly discussed above). Note that by Lemma \ref{lem:splitoflargevertices2}(ii), we have
$1 \leq \alpha \leq \delta - 1$.  It turns out that
the $\alpha = 1$ case is the most delicate, as that corresponds to an
extremal construction and cannot be resolved by our current approach.
The $\alpha = 2$ case will be easily
proved using Lemma \ref{lem:tau-by-rho} and \eqref{ineq:NEED:main}.  For $3 \leq \alpha \leq
\delta - 1$, observe that the left hand side of
\eqref{ineq:NEED:alpha-big} is convex in $\alpha$ in that range, while
the right hand side is linear in $\alpha$.  Therefore, once we verify
it for $\alpha = 3$ and $\alpha = \delta - 1$, all intermediate values
of $\alpha$ will follow by convexity.

\begin{lemma}
  \label{lem:alpha=2}
  Let $\alpha = 2$, and $\delta \ge 4$ be an even integer. Then
  \eqref{ineq:NEED:main} holds.
%  If that $\theta >
%  \frac{\delta}{\delta+1} n$, then $A_1 \cup A_2$ cannot contain exactly two
%  huge $\overline{A}$-degree vertices (or equivalently, $\alpha \neq 2$).
\end{lemma}

\begin{lemma}
  \label{lem:alpha=convex-extremes}
  Let $\delta \geq 6$ be an even integer, and let $n', \theta, \rho$ be
  non-negative real numbers satisfying
  \begin{displaymath}
    \theta \leq n' \,,
    \quad\quad
    \rho \leq n' - \theta \,.
  \end{displaymath}
  Then \eqref{ineq:NEED:alpha-big} holds when $\alpha = 3$, and when
  $\alpha = \delta - 1$.
\end{lemma}

Apart from the most delicate $\alpha = 1$ case, only the case
$(\delta, \alpha) = (4, 3)$ remains.  For both of these remaining
situations, we now take the refined approach alluded to earlier, where
we employ the entire degree sequence of $G[\overline{A}]$.  Since
$e(A) = o(m)$, the total number of edges satisfies
\begin{equation}
  \label{ineq:m-by-degrees}
  m
  \ge
  \alpha \theta + \rho + \frac{1}{2} (d_1 + 2d_2 + 3d_3 + \cdots) + o(m) \,,
\end{equation}
where $d_i$ is the number of vertices of degree $i$ in the induced
subgraph $G[\overline{A}]$.

Using this inequality and Lemma
\ref{lem:tau-by-degree-sequence}, we see that to prove \eqref{ineq:NEED:main}, it suffices to show
\begin{equation}
  \label{ineq:NEED:degrees}
  \theta
  + \frac{1}{2} \left(
    \frac{d_0}{1} + \frac{d_2}{3} + \frac{d_4}{5} + \cdots
  \right)
  \leq
  \frac{n'}{2}
  + \frac{1}{\delta + 1} \left(
    \alpha \theta + \rho + \frac{1}{2} \left(
      d_1 + 2d_2 + 3d_3 + \cdots
    \right)
  \right) \,,
\end{equation}
when the constraint \eqref{ineq:degrees-deficiency} on the
degree sequence holds.

\begin{lemma}
  \label{lem:optimization-alpha-1}
  Let $\delta$ be a positive integer, and let $n', \theta, \rho, d_0,
  d_1, \ldots$ be non-negative real numbers satisfying
  \begin{displaymath}
    \theta \leq n' \,,
    \quad\quad
    \rho \leq n' - \theta \,,
    \quad \text{and} \quad
    \sum_i d_i = n' \,,
  \end{displaymath}
  together with \eqref{ineq:degrees-deficiency}.  Then Inequality
  \eqref{ineq:NEED:degrees} holds for $(\delta, \alpha) = (4, 3)$, as
  well as for $\alpha = 1$ and arbitrary $\delta \geq 2$.
\end{lemma}

\noindent This final ingredient completes our proof of Theorem
\ref{thm:bisection-min}, modulo the intermediate lemmas.  \hfill
$\Box$

\subsection{Splitting large degree vertices}

\label{sec:splitoflargevertices}

In this section, we prove Lemma \ref{lem:splitoflargevertices2} by
establishing a stronger, but more technical result.  Let us
parameterize the balance between the $e(A_i, \overline{A})$ in the
optimal split by defining $\lambda$ such that $e(A_1, \overline{A}) =
\lambda m_1$ and $e(A_2, \overline{A}) = (1 - \lambda) m_1$.  Recall
that $m_1 = e(A, \overline{A})$, $m_2 = e(\overline{A})$, and note
that now
\begin{displaymath}
  \theta
  = e(A_1, \overline{A}) - e(A_2, \overline{A})
  = (2\lambda - 1)m_1 \,.
\end{displaymath}
By construction, we always have $\frac{1}{2} \leq \lambda \leq 1$.  As
it turns out, the particular location of $\lambda$ gives useful
information about the nature of the partition $A_1 \cup A_2$.

\begin{lemma}
  \label{lem:splitoflargevertices}
  Let $\kappa$ be an integer for which $\lambda >
  \frac{\kappa+1}{2\kappa+1}$.  Then both:
  \begin{description}
  \item[(i)] $|A_1| \leq \kappa$, and all vertices of $A_1$ have
    $\overline{A}$-degree at least $\theta$; and
  \item[(ii)] $A_2$ contains at most $\kappa-1$ vertices with
    $\overline{A}$-degree $\geq \theta$, and the
    $\overline{A}$-degrees of all remaining vertices sum to at most $n'
    - \theta$.
  \end{description}
\end{lemma}

\noindent \textbf{Proof.}\, Our argument is based on the fact that no
subset of $A$ has $\overline{A}$-degree sum strictly between $(1 -
\lambda)m_1$ and $\lambda m_1$.  Let $a = |A_1|$, and suppose for
contradiction that $a \geq \kappa+1$.  Then, by moving the vertex of lowest
$\overline{A}$-degree from $A_1$ to $A_2$, the remaining
$\overline{A}$-degree sum of $A_1$ would be under $\lambda m_1$, but at
least
\begin{displaymath}
  \frac{a-1}{a} \cdot \lambda m_1
  >
  \frac{\kappa}{\kappa+1} \cdot \frac{\kappa+1}{2\kappa+1} m_1
  =
  \frac{\kappa}{2\kappa+1} m_1
  >
  (1 - \lambda) m_1 \,,
\end{displaymath}
thereby contradicting the optimality of the $A_1 \cup A_2$ split.
Hence we indeed have $|A_1| \leq \kappa$.  Furthermore, if any vertex
of $A_1$ has $\overline{A}$-degree strictly less than $\theta$, then
we can also improve the split by moving it to $A_2$.  This establishes
(i).

For the first part of (ii), observe that if $A_2$ contains $\kappa$
vertices with $\overline{A}$-degree at least $\theta = (2 \lambda -
1)m_1$, then the $\overline{A}$-degree sum of $A_2$ is already at
least
\begin{displaymath}
  \kappa \cdot (2 \lambda - 1)m_1 > \frac{\kappa}{2\kappa+1}m_1 > (1-\lambda)m_1 \,,
\end{displaymath}
contradicting the definition of $\lambda$.

The final part of (ii) is slightly more involved.  Let $v$ be a vertex
of $A_1$ of minimal $\overline{A}$-degree, and let that degree be $d$.
We will show that the $\overline{A}$-degree sum of all vertices in
$A_2$ with $\overline{A}$-degree less than $\theta = (2\lambda -
1)m_1$ is at most $d - (2\lambda - 1)m_1$, which will suffice because
$d$ is clearly at most $|\overline{A}| = n'$.  Indeed, assume for
contradiction that this is not the case, and consider moving $v$ to
$A_2$.  The $\overline{A}$-degree sum of $A_1$ will fall to $\lambda
m_1 - d$, which must be at most $(1 - \lambda)m_1$ by the assumed
optimality of the $A_1 \cup A_2$ split.  Now consider moving the
vertices with $\overline{A}$-degree below $\theta$ from $A_2$ to
$A_1$, one by one.  If we moved all of them, then $A_1$ would achieve
a total $\overline{A}$-degree sum strictly greater than
\begin{displaymath}
  ( \lambda m_1 - d )
  +
  (d - (2\lambda - 1)m_1 )
  =
  (1 - \lambda) m_1 \,.
\end{displaymath}
Yet each vertex that we move has $\overline{A}$-degree strictly less
than $\theta = \lambda m_1 - (1 - \lambda) m_1$, so at some point in
the process, we must have a split where both $e(A_i, \overline{A})$
are strictly between $(1 - \lambda) m_1$ and $\lambda m_1$.  This
contradicts the optimality of the original $A_1 \cup A_2$ split, and
completes our proof.  \hfill $\Box$

\vspace{2mm}

We now observe that Lemma \ref{lem:splitoflargevertices2} follows
easily from Lemma \ref{lem:splitoflargevertices}.

\vspace{2mm}

\noindent \textbf{Proof of Lemma \ref{lem:splitoflargevertices2}.}\,
Since $\theta = (2\lambda - 1) m_1$, our assumption that $\theta >
\frac{m}{\delta + 1} \geq \frac{m_1}{\delta+1}$ gives $\lambda >
\frac{\delta + 2}{2(\delta + 1)}$.  Recalling that $\delta$ is even,
let $\kappa = \frac{\delta}{2}$, and observe that
$\frac{\kappa+1}{2\kappa + 1} = \frac{\delta + 2}{2(\delta+1)}$, so we
may apply Lemma \ref{lem:splitoflargevertices} with $\kappa =
\frac{\delta}{2}$, and Lemma \ref{lem:splitoflargevertices2} follows
immediately.  \hfill $\Box$

\subsection{Tight components and the degree sequence}

\label{sec:tau-by-degree-sequence}

Consider an arbitrary graph $H$.  For each $i$, let $d_i$ denote the
number of vertices in $H$ whose degree is exactly $i$.  In this
section, we prove Lemma \ref{lem:tau-by-degree-sequence}, which claims
that the number of tight components in $H$ is at most
\begin{displaymath}
  \frac{d_0}{1} + \frac{d_2}{3} + \frac{d_4}{5} + \cdots
\end{displaymath}

\noindent \textbf{Proof.}\, For each odd $i$, let $t_i$ denote the
number of tight components of $H$ which have exactly $i$ vertices.
Tight components have odd order and only contain even-degree vertices,
because for every vertex $v$ in a tight component $T$, there is a
perfect matching of $T \setminus \{v\}$, and $v$ has either 0 or 2
neighbors in each edge of the matching.  Hence tight components of
order $i$ are entirely composed of vertices of even degree $\leq i-1$.
We therefore have the inequalities:
\begin{align*}
  t_1 &\leq d_0 \\
  t_1 + 3t_3 &\leq d_0 + d_2 \\
  t_1 + 3t_3 + 5t_5 &\leq d_0 + d_2 + d_4 \\
  &\,\,\, \vdots
\end{align*}
We claim that even when the $t_i$ are allowed to take real values (as
opposed to integers), the maximum value of $t_1 + t_3 + \cdots$ is
achieved by the greedy algorithm which chooses $t_1 = \frac{d_0}{1}$,
$t_3 = \frac{d_2}{3}$, $t_5 = \frac{d_4}{5}$, etc.  This is the unique
solution when all of the inequalities are tight.

Indeed, assume for contradiction that we have an optimal real solution
in which the inequality $t_1 + 3t_3 + \cdots + i t_i \leq d_0 + d_2 +
\cdots + d_{i-1}$ is not tight.  Then, one can improve the solution by
increasing $t_i$ by some $(i+2) \varepsilon$ and decreasing $t_{i+2}$ by
$i \varepsilon$.  Since this perturbation keeps $t_1 + 3t_3 + \cdots +
(i+2) t_{i+2}$ constant, it is clearly still a feasible solution, but
it increases the objective $t_1 + t_3 + t_5 + \cdots$ by $2 \varepsilon >
0$, contradicting our assumed optimality.

Therefore, the number of tight components $t_1 + t_3 + t_5 + \cdots$
is at most $\frac{d_0}{1} + \frac{d_2}{3} + \frac{d_4}{5} + \cdots$,
as claimed.  \hfill $\Box$

\subsection{Tight components and the minimum degree condition}

\label{sec:tau-by-rho}

In this section, we prove Lemma \ref{lem:tau-by-rho}, which claims that
the number of tight components in $G[\overline{A}]$ is at most
$\frac{n' + \rho}{\delta - \alpha + 1}$.  Our main tool is Lemma
\ref{lem:tau-by-degree-sequence}, which implies the following weaker
bound on $\tau$ in terms of the degree sequence of $G[\overline{A}]$:
\begin{displaymath}
  \tau \leq \frac{d_0}{1} + \frac{d_1}{2} + \frac{d_2}{3} + \cdots
\end{displaymath}
(Recall that $d_i$ is the number of vertices of $G[\overline{A}]$ that
have degree exactly $i$.)

To get an effective bound on $\tau$, we would like to see how large the
right hand side of the above can be.
Our constraints are \eqref{ineq:degrees-deficiency}, which we restate here,

\begin{equation}
\label{ineq:degrees-deficiency-re}
  (\delta - \alpha) d_0
  + (\delta - \alpha - 1) d_1
  + \cdots
  + d_{\delta - \alpha - 1}
  \leq \rho \, ,
\end{equation}
together with the
following obvious constraint (which we relax from equality)
\begin{equation}
  \label{ineq:n'-vtxs}
  \sum_i d_i \leq n' \,.
\end{equation}
We claim that when the $d_i$ are allowed to take any non-negative
real values (as opposed to integers), the maximum value of
$\frac{d_0}{1} + \frac{d_1}{2} + \cdots$ subject to the constraints
\eqref{ineq:degrees-deficiency-re} and \eqref{ineq:n'-vtxs} is achieved
when $d_0$ and $d_{\delta-\alpha}$ are the only nonzero $d_i$.

To see this, consider an optimizer, and note that if any $i > \delta - \alpha$
has $d_i > 0$, then we can maintain the left hand sides of
\eqref{ineq:degrees-deficiency-re} and \eqref{ineq:n'-vtxs} by reducing
$d_i$ to 0, and increasing $d_{\delta-\alpha}$ by the same amount.
This increases the value of the objective.  Next observe that if any
$0 < i < \delta-\alpha$ has $d_i > 0$, then we can maintain the left
hand side of \eqref{ineq:degrees-deficiency-re} by reducing $d_i$ to 0,
and increasing $d_0$ by $\frac{\delta-\alpha - i}{\delta-\alpha} d_i$.
This only reduces the left hand side of \eqref{ineq:n'-vtxs}, but
increases the objective by $\frac{\delta-\alpha - i}{\delta-\alpha}
d_i - \frac{d_i}{i+1}$.  Since $\frac{\delta-\alpha-i}{\delta-\alpha}
\geq \frac{1}{i+1}$ for all $0 \leq i \leq \delta - \alpha - 1$, we see that
this cannot decrease the objective.

Therefore, we may assume that our optimizer is only supported at $d_0$
and $d_{\delta-\alpha}$.  Inequality \eqref{ineq:degrees-deficiency}
then simplifies to $d_0 \leq \frac{\rho}{\delta-\alpha}$.  In the
objective, the coefficient of $d_0$ is strictly greater than that of
$d_\delta$, hence to maximize the objective, we must choose $d_0$ to be as
large as possible, i.e., $d_0 = \frac{\rho}{\delta-\alpha}$, and then
choose maximum $d_{\delta-\alpha}$ subject to the constraint \eqref{ineq:n'-vtxs},
which says $d_0 + d_{\delta-\alpha} \le n'$.
Therefore, the objective is at most
\begin{displaymath}
  \frac{\rho/(\delta-\alpha)}{1}
  + \frac{n' - (\rho/(\delta-\alpha))}{\delta-\alpha + 1}
  =
  \frac{n' + \rho}{\delta-\alpha + 1} \,,
\end{displaymath}
as claimed.  \hfill $\Box$

\subsection{Two huge degree vertices}

In this section, we prove Lemma \ref{lem:alpha=2}, which considers
the case when $A$ contains exactly two huge degree vertices (vertices with
$\overline{A}$-degree at least $\theta$).

\vspace{2mm}

%\begin{equation}
%  \label{ineq:NEED:main}
%  \theta + \frac{\tau}{2} \leq \frac{n}{2} + \frac{m}{\delta + 1}
%  \,.
%\end{equation}

\noindent \textbf{Proof of Lemma \ref{lem:alpha=2}.}\, We will first
prove that $\theta \le \frac{3}{5}n'$. Suppose not.
Lemma \ref{lem:splitoflargevertices2}(iii) says that the sum of
$\overline{A}$-degrees of the non-huge degree vertices in $A_2$ is
$\rho \leq n' - \theta < \theta$. Since, by construction, $e(A_1, \overline{A})
\geq e(A_2, \overline{A})$, it is impossible for
both of the huge $\overline{A}$-degree vertices to be in $A_1$, because it implies that
\begin{displaymath}
  \theta = e(A_1, \overline{A}) - e(A_2, \overline{A})
  \ge
  2 \theta - \rho
  > 
  \theta \,.
\end{displaymath}
which is a contradiction. Thus, the only remaining case is when $A_1$ and $A_2$ each have a
huge-degree vertex.  Then since $e(A_1, \overline{A}) \le n'$
and $e(A_2, \overline{A}) \ge \theta$, the maximum possible gap is
\begin{displaymath}
  \theta = e(A_1, \overline{A}) - e(A_2, \overline{A})
  \leq
  n' - \theta < \theta \,,
\end{displaymath}
which again is a contradiction.

Now, using $\theta \le \frac{3}{5}n'$ and Lemma \ref{lem:tau-by-rho} with
$\alpha = 2$, $\delta \ge 4$,
we will verify \eqref{ineq:NEED:main} which says,
\[ \theta + \frac{\tau}{2} \leq \frac{n}{2} + \frac{m}{\delta + 1}. \]
Note that by Lemmas \ref{lem:splitoflargevertices2} and \ref{lem:tau-by-rho}, we have
$\tau \le \frac{n' + \rho}{\delta -\alpha+ 1}=
\frac{n' + \rho}{\delta - 1} \le \frac{2n'-\theta}{\delta - 1} \le \frac{2n' - \theta}{3}$.
By the minimum degree assumption, we also have $m \ge \frac{\delta}{2}n'$.
Thus it suffices to prove that
\[\theta + \frac{2n' - \theta}{6}
  \leq
  \frac{n'}{2}
  + \frac{\delta n'}{2(\delta + 1)}, \]
and since $\delta \ge 4$, that
$\frac{5}{6}\theta
  \leq
  \left( \frac{1}{2}  + \frac{2}{5} - \frac{1}{3} \right)n' = \frac{17}{30}n'$.
This is true as $\theta \le \frac{3}{5}n'$. \hfill $\Box$

\subsection{Convexity extremes}

In this section, we prove Lemma \ref{lem:alpha=convex-extremes}, which
disposes of the $\alpha = 3$ and $\alpha = \delta - 1$ cases of
\eqref{ineq:NEED:alpha-big} when $\delta \geq 6$.  We begin with
$\alpha = 3$, which corresponds to

\begin{lemma}
  Let $\delta \geq 5$, and let $n', \theta, \rho$ be non-negative real
  numbers satisfying $\theta \leq n'$ and $\rho \leq n' - \theta$.
  Then
  \begin{displaymath}
    \theta + \frac{n' + \rho}{2(\delta - 2)}
    \leq
    \frac{2\delta + 1}{2(\delta+1)} n'
    + \frac{3\theta + \rho}{2(\delta+1)} \,.
  \end{displaymath}
\end{lemma}

\noindent \textbf{Proof.}\, The coefficient of $\rho$ on the left hand
side is clearly larger than its coefficient on the right hand side, so
the inequality is sharpest when $\rho$ hits its constraint of $n' -
\theta$.  It then suffices to show
\begin{displaymath}
  \theta + \frac{2n' - \theta}{2(\delta - 2)}
  \leq
  \frac{2\delta + 1}{2(\delta+1)} n'
  + \frac{2\theta + n'}{2(\delta+1)} \,.
\end{displaymath}
Next, observe that the coefficient of $\theta$ on the left hand side
is $\frac{2\delta-5}{2(\delta-2)}$, while its coefficient on the right
hand side is $\frac{2}{2(\delta+1)}$.  The left-hand coefficient is
easily larger in our range $\delta \geq 5$, so the inequality sharpens
when we increase $\theta$ to its constraint of $n'$.  It therefore
remains to show that
\begin{displaymath}
  n' + \frac{n'}{2(\delta-2)}
  \leq
  \frac{2\delta + 1}{2(\delta+1)} n'
  + \frac{3n'}{2(\delta+1)} \,,
\end{displaymath}
or equivalently,
\begin{displaymath}
  \frac{2\delta-3}{2(\delta-2)}
  \leq
  \frac{2\delta+4}{2(\delta+1)} \,,
\end{displaymath}
which one easily verifies is true for all $\delta \geq 5$.
\hfill $\Box$

\vspace{2mm}

Next we establish the $\alpha = \delta-1$ case.

\begin{lemma}
  Let $\delta \geq 5$, and let $n', \theta, \rho$ be non-negative real
  numbers satisfying $\theta \leq n'$ and $\rho \leq n' - \theta$.
  Then
  \begin{displaymath}
    \theta + \frac{n' + \rho}{4}
    \leq
    \frac{2\delta + 1}{2(\delta+1)} n'
    + \frac{(\delta-1)\theta + \rho}{2(\delta+1)} \,.
  \end{displaymath}
\end{lemma}

\noindent \textbf{Proof.}\, The coefficient of $\rho$ on the left hand
side is clearly greater than its coefficient on the right hand side, so
we may again assume that $\rho = n' - \theta$.  It remains to show
\begin{displaymath}
  \theta + \frac{2n' - \theta}{4}
  \leq
  \frac{2\delta + 1}{2(\delta+1)} n'
  + \frac{(\delta-2)\theta + n'}{2(\delta+1)} \,.
\end{displaymath}
The coefficient of $\theta$ on the left hand side is $\frac{3}{4}$,
which is again easily greater than its coefficient of
$\frac{\delta-2}{2(\delta+1)}$ on the right hand side; hence we may
assume $\theta = n'$, leaving us to establish
\begin{displaymath}
  n' + \frac{n'}{4}
  \leq
  \frac{2\delta + 1}{2(\delta+1)} n'
  + \frac{(\delta-1)n'}{2(\delta+1)} \,,
\end{displaymath}
which is equivalent to
\begin{displaymath}
  \frac{5}{4}
  \leq
  \frac{3\delta}{2(\delta+1)} \,.
\end{displaymath}
One easily sees that this indeed holds for $\delta \geq 5$.  \hfill
$\Box$

\subsection{Tight components with three huge degree vertices}

\label{sec:alpha=3}

We prove Lemma \ref{lem:optimization-alpha-1} for the case $\delta =
4$ and $\alpha = 3$.  Using these parameters in
\eqref{ineq:NEED:degrees} and \eqref{ineq:degrees-deficiency}, we see
that we need to establish
\begin{equation}
  \label{ineq:alpha=3}
  \theta
  + \frac{1}{2} \left(
    \frac{d_0}{1} + \frac{d_2}{3} + \frac{d_4}{5} + \cdots
  \right)
  \leq
  \frac{n'}{2}
  + \frac{1}{5} \left(
    3 \theta + \rho + \frac{1}{2} \left(
      d_1 + 2d_2 + 3d_3 + \cdots
    \right)
  \right) \,,
\end{equation}
subject to the constraints
\begin{equation}
  \label{ineq:alpha=3-constraint}
  d_0 \leq \rho \,,
  \quad\quad
  \theta \leq n' \,,
  \quad\quad
  \rho \leq n' - \theta \,,
  \quad\text{and}\quad
  \sum_i d_i = n' \,.
\end{equation}
We consider the differences between the coefficients of the $d_i$ on
the left and on the right.  For $d_0$, its coefficient on the left is
$\frac{1}{2}$, while its coefficient on the right is 0.  This is the
only coefficient where that differential is in favor of the left hand
side.  Indeed, it is clear that for even $i$, the difference between
the left hand coefficient and the right hand coefficient strictly
decreases as $i$ increases.  A similar fact is true for odd $i$.  Note
that the differential for the coefficient of $d_1$ is $0 -
\frac{1}{10}$, whereas the differential for $d_2$ is $\frac{1}{6} -
\frac{1}{5} = -\frac{1}{30}$.  Hence the differential is largest
for $d_2$ (apart from $d_0$) since $-\frac{1}{30} > -\frac{1}{10}$.

Therefore, in view of the constraints $d_0 \leq \rho$ and $\sum_i d_i
= n'$, inequality \eqref{ineq:alpha=3} is sharpest when $d_0$ is as
large as possible, and all the rest of the weight is on $d_2$.  It
remains to show that
\begin{displaymath}
  \theta + \frac{1}{2} \left( \rho + \frac{n' - \rho}{3} \right)
  \leq
  \frac{n'}{2} + \frac{1}{5} \left(
    3\theta + \rho + (n' - \rho)
  \right) \,,
\end{displaymath}
or equivalently,
\begin{displaymath}
  \theta + \frac{\rho}{3} + \frac{n'}{6}
  \leq
  \frac{7}{10} n' + \frac{3}{5} \theta \,.
\end{displaymath}
Since $\rho$ only appears on the left hand side, the inequality is
sharpest when $\rho$ hits its upper constraint of $n'-\theta$; it
remains to show
\begin{displaymath}
  \frac{2}{3} \theta + \frac{n'}{2}
  \leq
  \frac{7}{10} n' + \frac{3}{5} \theta \,.
\end{displaymath}
This is equivalent to $\frac{1}{15}\theta \le \frac{1}{5}n'$, which
clearly holds since $\theta \le n'$.  \hfill $\Box$

\subsection{Tight components with a single vertex of huge degree}

We now prove our sharpest estimate (Lemma
\ref{lem:optimization-alpha-1}), which handles the situation when $A$
contains a single vertex of huge degree.  As we need to establish
\eqref{ineq:NEED:degrees}, we rearrange it to isolate the involvement
of the degree sequence.  We must show:
\begin{equation}
  \label{ineq:NEED:degrees-on-left}
  \frac{1}{2} \left(
    \frac{d_0}{1} + \frac{d_2}{3} + \frac{d_4}{5} + \cdots
  \right)
  - \frac{1}{2(\delta+1)} \left(
    d_1 + 2d_2 + 3d_3 + \cdots
  \right)
  \leq
  \frac{n'}{2} + \frac{\theta + \rho}{\delta + 1} - \theta \,.
\end{equation}

We first control the left hand side in terms of $\rho$, using the
given constraint \eqref{ineq:degrees-deficiency}.

\begin{lemma}
  \label{lem:opt-degrees-rho}
  For every positive even integer $\delta$ and non-negative real
  numbers $n', \rho, d_0, d_1, \ldots$, the maximum value of the
  following objective:
  \begin{displaymath}
    \frac{1}{2} \left(
      \frac{d_0}{1} + \frac{d_2}{3} + \frac{d_4}{5} + \cdots
    \right)
    - \frac{1}{2(\delta+1)} \left(
      d_1 + 2d_2 + 3d_3 + \cdots
    \right) \,,
  \end{displaymath}
  subject to the constraints
  \begin{displaymath}
    \sum_{i=0}^{\delta-2} (\delta - 1 - i) d_i \leq \rho \,,
    \quad
    \text{and}
    \quad
    \sum_i d_i = n' \,,
  \end{displaymath}
  is at most
  \begin{equation}
    \label{ineq:obj-by-rho}
    - \frac{\delta-1}{2(\delta+1)} n' + \frac{\delta}{(\delta - 1) (\delta + 1)} \rho \,.
  \end{equation}
\end{lemma}

\noindent \textbf{Proof.}\, Let us study the coefficients of the $d_i$
in the objective.  Coincidentally, the coefficients of $d_{\delta-1}$
and $d_\delta$ are both equal to
\begin{displaymath}
  0 - \frac{1}{2(\delta+1)} (\delta - 1)
  =
  - \frac{\delta - 1}{2(\delta+1)}
  =
  \frac{1}{2} \cdot \frac{1}{\delta+1} - \frac{1}{2(\delta+1)} \cdot \delta \,,
\end{displaymath}
while all coefficients of $d_i$ with $i \geq \delta + 1$ are strictly
smaller.  Therefore, by moving weight from $d_i$ to $d_{\delta-1}$ so that
$\sum_i d_i$ remains constant, we can increase the objective function without
violating any constraints. This allows us to conclude that there is an
optimal point which has $d_i = 0$ for all $i \geq \delta$.

We aim to show that there is in fact an optimizer which is supported
only at $d_0$ and $d_{\delta-1}$.  So, suppose that the optimizer has
$d_i > 0$ for some $1 \leq i \leq \delta - 2$.  Consider reducing
$d_i$ to 0, increasing $d_0$ by $\frac{\delta - 1 - i}{\delta - 1}
d_i$, and increasing $d_{\delta-1}$ by $\frac{i}{\delta-1} d_i$.
(This preserves both $\sum_{i=0}^{\delta - 2} (\delta - 1 - i) d_i$
and $\sum_i d_i$.)  Note that the coefficients of $d_0$ and
$d_{\delta-1}$ in the objective are $\frac{1}{2}$ and
$-\frac{\delta-1}{2(\delta+1)}$, respectively.

\vspace{2mm}

\noindent \textbf{Case 1: $\boldsymbol{i}$ is odd.}\, The coefficient
of $d_i$ in the objective is $-\frac{i}{2(\delta+1)}$, so this
perturbation increases the objective by:
\begin{align*}
  & d_i \left[
    \frac{1}{2} \cdot \frac{\delta - 1 - i}{\delta - 1}
    - \frac{i}{2(\delta+1)} \cdot (-1)
    + \left( -\frac{\delta-1}{2(\delta+1)} \right) \cdot \frac{i}{\delta-1}
  \right] \\
  & \quad\quad\quad = d_i \left[
    \frac{
      (\delta-1-i)(\delta+1) + i(\delta-1) - (\delta-1)i
    }
    {2(\delta-1)(\delta+1)}
  \right] \\
  & \quad\quad\quad = d_i \cdot \frac{\delta-1-i}{2(\delta-1)} \,,
\end{align*}
which is clearly positive for $i \leq \delta-2$.

\vspace{2mm}

\noindent \textbf{Case 2: $\boldsymbol{i}$ is even.}\, The coefficient
of $d_i$ in the objective is now $\frac{1}{2(i+1)} -
\frac{i}{2(\delta+1)}$, so by the previous calculation this
perturbation increases the objective by:
\begin{displaymath}
  d_i \cdot \frac{\delta-1-i}{2(\delta-1)} - d_i \cdot \frac{1}{2(i+1)} \,,
\end{displaymath}
which is non-negative precisely when
\begin{equation}
  \label{ineq:convexity-surprise}
  \frac{\delta-1-i}{2(\delta-1)} \geq \frac{1}{2(i+1)} \,.
\end{equation}
Yet it is clear that \eqref{ineq:convexity-surprise} has equality at
both $i = 0$ and $i = \delta - 2$.  Since the left hand side is
linear in $i$, while the right hand side is convex on the domain $i >
-1$, we immediately conclude that \eqref{ineq:convexity-surprise}
holds for all $0 \leq i \leq \delta - 2$.

\vspace{2mm}

As both the odd and even cases produce perturbations that do not
decrease the objective, we may therefore consider an optimizer which
is entirely supported at $d_0$ and $d_{\delta - 1}$.  The constraints
then translate into
\begin{displaymath}
  (\delta - 1) d_0 \leq \rho
  \quad \text{and} \quad
  d_0 + d_{\delta-1} = n' \,.
\end{displaymath}
Since the coefficient of $d_0$ in the objective $\big( \frac{1}{2}
\big)$ is strictly larger than that of $d_{\delta-1}$ (which is $-
\frac{\delta - 1}{2(\delta + 1)})$, to maximize the objective function,
we need to take $d_0$ to be as large as possible, i.e., $d_0 = \frac{\rho}{\delta-1}$,
and then choose $d_{\delta-1}$ to satisfy $d_0 + d_{\delta-1} = n'$.
Hence we conclude that the objective is indeed bounded from above by
\begin{displaymath}
  \frac{1}{2} \cdot \frac{\rho}{\delta - 1}
  - \frac{\delta - 1}{2(\delta + 1)} \cdot \left(
    n' - \frac{\rho}{\delta - 1}
  \right) \,,
\end{displaymath}
which is precisely \eqref{ineq:obj-by-rho}.   \hfill $\Box$

\vspace{2mm}

Our final ingredient for the proof of Theorem \ref{thm:bisection-min}
now follows easily.

\vspace{2mm}

\noindent \textbf{Proof of Lemma \ref{lem:optimization-alpha-1}.}\,
Section \ref{sec:alpha=3} establishes the result for the case $(\delta,
\alpha) = (4, 3)$, so it remains to consider the case when $\alpha =
1$.  Using the result of Lemma \ref{lem:opt-degrees-rho} in
\eqref{ineq:NEED:degrees-on-left}, we see that it remains to establish
\begin{displaymath}
  - \frac{\delta-1}{2(\delta+1)} n' + \frac{\delta}{(\delta - 1) (\delta + 1)} \rho
  \leq
  \frac{n'}{2} + \frac{\theta + \rho}{\delta + 1} - \theta \,,
\end{displaymath}
over the domain bounded by $0 \leq \theta \leq n'$ and $\rho \leq n' -
\theta$.  Equivalently, we need
\begin{displaymath}
  \frac{\delta}{\delta + 1} \theta
  +
  \frac{1}{(\delta-1)(\delta+1)} \rho
  \leq
  \frac{\delta}{\delta + 1} n' \,.
\end{displaymath}
Using $\rho \leq n' - \theta$ to eliminate $\rho$, we see that it
suffices to prove
\begin{displaymath}
  \left(
    \frac{\delta}{\delta + 1} - \frac{1}{(\delta-1)(\delta+1)}
  \right) \theta
  \leq
  \left(
    \frac{\delta}{\delta+1} - \frac{1}{(\delta-1)(\delta+1)}
  \right) n' \,,
\end{displaymath}
which clearly follows from $\theta \leq n'$.  This completes our
proof.  \hfill $\Box$

\section{Analysis of extremal examples}
\label{sec:bisection-min-tight}

In this section, we establish the asymptotic extremality of the two
families of constructions which we described after Theorem
\ref{thm:bisection-min} in the Introduction.

\vspace{2mm}

\noindent \textbf{Proof of Proposition
  \ref{prop:bisection-min-tight}.}\, Let $\delta \geq 2$ be a fixed
even number.  We begin with the first family.  Let $x$ and $y$ be
non-negative integers, with $y$ odd.  Let $G = (V, E)$ be the
vertex-disjoint union of $x$ copies of $K_\delta$ and $y$ copies of
$K_{\delta + 1}$, together with a new vertex $v_0$ which is adjacent
to all other vertices.  The number of vertices is then
\begin{displaymath}
  n = \delta x + (\delta+1) y + 1 \,,
\end{displaymath}
which is even, and the number of edges is
\begin{displaymath}
  m
  = \binom{\delta}{2} x + \binom{\delta+1}{2} y + \delta x + (\delta+1)y
  = \frac{\delta (\delta+1)}{2} x + \frac{(\delta+1)(\delta+2)}{2} y \,.
\end{displaymath}
Let $V = V_1 \cup V_2$ be an arbitrary bisection, and without loss of
generality, suppose that the dominating vertex $v_0$ lies in $V_1$.
The induced subgraph $G[V_1 \setminus \{v_0\}]$ is then the union of
$t \leq x+y$ cliques, with orders $c_1, c_2, \ldots, c_t$.  Since we
have a bisection, $\sum c_i = \frac{n}{2} - 1$.  Using convexity, it
is then clear that $\sum \binom{c_i}{2}$ is minimized when $t = x+y$,
and the $c_i$ are within 1 of each other. Since $\frac{n}{2} -1 =
\frac{\delta}{2}(x+y) + \frac{y-1}{2}$, this happens precisely when
$x + \frac{y+1}{2}$ of the $c_i$'s are $\frac{\delta}{2}$, and
$\frac{y-1}{2}$ of the $c_i$'s are $\frac{\delta}{2} + 1$.  Therefore,
\begin{eqnarray*}
  e(V_1)
  &=&
  \left( \frac{n}{2} - 1 \right) + \sum_{i=1}^t \binom{c_i}{2} \\
  &\geq&
  \left( \frac{\delta x + (\delta+1)y + 1}{2} - 1 \right)
  + \left( x + \frac{y+1}{2} \right) \binom{\delta/2}{2}
  + \frac{y-1}{2} \binom{(\delta/2) + 1}{2} \\
  &=&
  x \left[ \frac{\delta}{2} + \binom{\delta/2}{2} \right]
  + y \left[
    \frac{\delta+1}{2} + \frac{1}{2} \binom{\delta/2}{2}
    + \frac{1}{2} \binom{(\delta/2) + 1}{2}
  \right] \\
  & &  + \left[
    -\frac{1}{2} + \frac{1}{2} \binom{\delta/2}{2}
    - \frac{1}{2} \binom{(\delta/2)+1}{2}
  \right] \\
  &=&
  x \cdot \frac{ \delta (\delta + 2) }{8}
  + y \cdot \frac{ (\delta + 2)^2 }{8}
  - \frac{\delta + 2}{4} \\
  &=&
  \frac{\delta+2}{4(\delta+1)} m - \frac{\delta + 2}{4} \,,
\end{eqnarray*}
which is indeed $\big( \frac{\delta+2}{4(\delta+1)} - o(1) \big) m$
since $\delta$ is constant.

Proceeding to the second family, let $G = (V, E)$ be the complete
bipartite graph $K_{\delta+1, n-\delta-1}$.  The number of edges is
\begin{displaymath}
  m = (\delta + 1) (n - \delta - 1) \,.
\end{displaymath}
Consider a bisection $V = V_1 \cup V_2$.  Of the $\delta+1$ vertices
from the smaller part of the $K_{\delta+1, n-\delta-1}$, suppose that
$x$ of them are in $V_1$ and $y$ of them are in $V_2$.  Without loss
of generality, suppose that $x \geq y$; then, $x \geq \frac{\delta}{2}
+ 1$, and
\begin{displaymath}
  e(V_1)
  =
  x \left( \frac{n}{2} - x \right)
  \geq
  \left( \frac{\delta}{2} + 1 \right)
  \left( \frac{n}{2} - \frac{\delta}{2} - 1 \right)
  =
  \frac{\delta+2}{4(\delta+1)}m - \frac{\delta+2}{4} \,.
\end{displaymath}
Note that the bound we obtained from the two families are exactly the same.

%The ratio is asymptotically
%\begin{displaymath}
%  \frac{e(V_1)}{m}
%  \geq
%  \frac{\left( \frac{\delta}{2} + 1 \right)
%    \left( \frac{n}{2} - \frac{\delta}{2} - 1 \right)}
%  { (\delta + 1) (n - \delta - 1) }
%  =
%  (1 - o(1)) \frac{\delta+2}{4(\delta+1)} \,,
%\end{displaymath}
%as claimed.  \hfill $\Box$

\section{Related bisection results}
\label{sec:related-proofs}

In this section, we prove the results which were stated in Section
\ref{sec:related-results}.  We begin by proving that $K_{n/3, 2n/3}$
has a bisection of size at least $\frac{3}{4} m$, but every bisection
induces at least $\frac{m}{4}$ edges on some side.

\vspace{2mm}

\noindent \textbf{Proof of Proposition
  \ref{prop:big-bisection-no-judicious}.}\, In the graph
$K_{n/3,2n/3}$, let $A$ be the part of size $\frac{n}{3}$ and $B$ be
the part of size $\frac{2n}{3}$. The total number of edges is
$m=\frac{2n^{2}}{9}$. The largest bisection is obtained by placing all
the vertices of $A$ on one side and arbitrarily allocating the
vertices of $B$ to create a bisection.  This produces a bisection of
size $\frac{n}{3}\cdot\frac{n}{2}=\frac{n^{2}}{6}=\frac{3}{4}m$.

Now consider a bisection $V_{1}\cup V_{2}$. Without loss of generality
we may assume that $V_{1}$ contains at least $xn \ge \frac{|A|}{2} =
\frac{n}{6}$ vertices of $A$ (thus $\frac{1}{6}\le
x\le\frac{1}{3}$). Then the number of edges in $V_{1}$ is
\begin{displaymath}
  e(V_1)
  =
  xn \cdot \left(
    \frac{n}{2} - xn
  \right)
  =
  x \left(\frac{1}{2}-x\right) n^2
  \geq
  \frac{1}{18} n^2
  =
  \frac{m}{4} \,.
\end{displaymath}
Therefore, $\max\{ e(V_{1}),e(V_{2})\} \geq \frac{m}{4}$ for every
bisection of the vertex set.  \hfill $\Box$

\vspace{2mm}

Next, we prove Theorem \ref{thm:max-degree-bisection}, which finds
almost-bisections of size at least $\frac{r+1}{2r} m$ or
$\frac{r+2}{2(r+1)} m$ (depending on the parity of $r$) in graphs with
degrees bounded by $r$.

\vspace{2mm}

\noindent \textbf{Proof of Theorem \ref{thm:max-degree-bisection}.}\,
Let $G = (V, E)$ be a graph of maximum degree at most $r$.  The
Hajnal-Szemer\'edi theorem \cite{HaSz70}, produces an \emph{equitable}\/ coloring of
$G$ with $r+1$ colors, i.e., a coloring where all color class sizes differ by at
most one. Let $W_1, W_2, \ldots, W_{r+1}$ be the color classes.

Consider first the case when $r$ is odd.  Generate a random
bipartition $V = V_1 \cup V_2$ by uniformly selecting exactly
$\frac{r+1}{2}$ of the color classes to put in $V_1$, and placing the
other $\frac{r+1}{2}$ color classes in $V_2$.  The probability that a
fixed edge contributes towards $e(V_{1},V_{2})$ is exactly
$(\frac{r+1}{2})^{2}/{r+1 \choose 2}=\frac{r+1}{2r}$, so by linearity
of expectation, the expected value of $e(V_{1,}V_{2})$ is
$\frac{r+1}{2r}m$.  Consequently, there exists a bisection which
contains at least so many edges.  Since all color classes differ in
size by at most 1, we automatically have that the $V_i$ differ in size
by at most $\frac{r+1}{2}$, as claimed.

Now consider the case when $r$ is even. Generate a random bisection
$V_1 \cup V_2$ of $G$ as follows.  First, pick a special index $k \in
[r+1]$ uniformly at random, generate a random bisection of $W_k$, and
put the parts into $V_1$ and $V_2$. Then, uniformly distribute the
remaining color classes as in the odd case above.  There are two ways
for a fixed edge $e$ to contribute to $e(V_1, V_2)$.  The probability
that one endpoint lies in $W_k$ is $\frac{2}{r+1}$, and then the other
endpoint has probability $\frac{1}{2}$ of being distributed onto the
other side of the bipartition.  On the other hand, the probability
that neither endpoint lies in $W_k$, and the endpoint color classes
are separated across the bipartition, is exactly
$\frac{r-1}{r+1} \cdot \frac{r^2/4}{{r \choose 2}} = \frac{r}{2(r+1)}$.
Therefore, the probability that $e$ contributes to
$e(V_1, V_2)$ is exactly
\begin{displaymath}
  \frac{2}{r+1} \cdot \frac{1}{2}
  +
  \frac{r}{2(r+1)}
  =
  \frac{r+2}{2(r+1)} \,.
\end{displaymath}
Consequently, the expected value of $e(V_{1},V_{2})$ is
$\frac{r+2}{2(r+1)}m$, and there exists a bisection which contains at
least so many edges.  Again, the gap between the $|V_i|$ is at most $1
+ \frac{r}{2}$, because the bisection of $W_k$ can introduce an error
of at most 1, and the other $r$ color classes differ in size by at
most 1.  \hfill $\Box$

\vspace{2mm}

The previous bound easily produces a result on perfect bisections of
bounded degree graphs.

\begin{corollary}
If $G$ has maximum degree at most $r$, then there exists a bisection
of size at least $\frac{r+1}{2r}m - \frac{r(r+1)}{4}$ if $r$ is
odd and $\frac{r+2}{2(r+1)}m - \frac{r(r+2)}{4}$ if $r$ is even.
\end{corollary}

\noindent \textbf{Proof.}\, If $r$ is odd, the previous theorem
produces a bipartition of size at least $\frac{r+1}{2r} m$, and sides
differing in size by at most $\frac{r+1}{2}$.  We may then balance the
bipartition by moving at most $\frac{r+1}{4}$ vertices.  Since all
degrees are bounded by $r$, this affects the size of the cut by at
most $\frac{r(r+1)}{4}$, as claimed.

On the other hand, if $r$ is even, then the previous theorem produces
a bipartition of size at least $\frac{r+2}{2(r+1)} m$, which can be
balanced by moving at most $\frac{r+2}{4}$ vertices.  By a similar
argument to above, this affects the cut size by at most
$\frac{r(r+2)}{4}$, as claimed.  \hfill $\Box$

\vspace{2mm}

Finally, we prove that every $r$-regular graph has a bisection of size at least
$\frac{r+1}{2r}m$ when $r$ is odd and $\frac{r+2}{2(r+1)}m$ when $r$ is
even.

\vspace{2mm}

\noindent \textbf{Proof of Theorem \ref{thm:regular-bisection}.}\,
An $r$-regular graph clearly has chromatic number $k \le r+1$. Thus
by the relation between Max Cut and the chromatic number of a graph
mentioned in the introduction (see the discussion following Proposition
\ref{prop:big-bisection-no-judicious}, or \cite{AnGrLi83, LeTu82, Locke82}),
we can find a bipartition $V = V_1 \cup V_2$ with the required size,
but where the part sizes are not necessarily the same.
It suffices to prove that we can always move
vertices from the larger part ($V_1$, say) to the smaller part ($V_2$,
say) while keeping the property that the bipartition has at least the
required size.  Indeed, if there is a vertex $v \in V_1$ with at least
as many neighbors in $V_1$ as in $V_2$, then moving $v$ to $V_2$
would not decrease the size of the bipartition. Repeatedly move such vertices until
either we reach a bisection and are done, or there
are no more such vertices.

In the later case, every vertex in $V_1$
has more neighbors in $V_2$ than in $V_1$. Move $|V_1| - \frac{n}{2}$
vertices of $V_1$ to $V_2$. Note that the number of neighbors in $V_2$ of
any of the remaining $\frac{n}{2}$ vertices of $V_1$ only grows. Therefore,
each of them still have strictly more than $\frac{r}{2}$ neighbors in $V_2$.
Yet $G$ was $r$-regular, so $m = \frac{rn}{2}$.  If $r$ is odd, the
number of crossing edges is then at least $\frac{n}{2} \cdot
\frac{r+1}{2} = \frac{r+1}{2r} m$.  If $r$ is even, the number of
crossing edges is at least $\frac{n}{2} \cdot \big( \frac{r}{2} + 1
\big) = \frac{r+2}{2r}m > \frac{r+2}{2(r+1)}m$.  Therefore, in both cases we
obtain a bisection, while still maintaining
the desired cut size.  \hfill $\Box$

\section{Concluding remarks}

In this paper, we studied the graph bisection problem from several
different angles.  Theorem~\ref{thm:tight} extended the classical
Edwards bound for graph bipartition by introducing the maximum degree,
together with a new parameter, the number of tight components.  These
tight components were very useful in our analysis, but have not been
widely studied in the literature.  Indeed, we are not aware of results
concerning their general structure.  For example, the simplest example
of a tight component is a clique of odd order, and more complex tight
components can be assembled by taking two vertex-disjoint tight
components, and identifying one vertex from each of them.  Yet it is
not clear if every tight component can be obtained through such a
process.

Using our parameterization in terms of the number of tight components,
we then pushed our approach further to study judicious bisections of
graphs with respect to minimum-degree constraints
(Theorem~\ref{thm:bisection-min}).  Part of the difficulty in proving
this result stemmed from the wide range of extremal examples, as
described in Proposition~\ref{prop:bisection-min-tight}.  From our
proof, it appears that these examples are in some sense exhaustive.
More precisely, when the inequalities are tight (up to the $o(m)$
error terms), we either have (i) exactly $\delta+1$ large degree
vertices, each of degree about $n$, with no other edges, or (ii)
exactly one vertex of large degree, with degree about $n$, and with
the correct number of tight components in the remainder.  However,
since the error term still remains in Theorem~\ref{thm:bisection-min},
we leave this as only an observation, and ask whether
one can remove the error term from Theorem~\ref{thm:bisection-min}.
More precisely, is it true that for every positive even integer
$\delta$, every graph of minimum degree at least $\delta$ admits a
bisection in which the number of edges within each part is at most
$\frac{\delta + 2}{4(\delta+1)}m$?

For directed graphs, one can also ask questions in which one seeks to
maximize or minimize several quantities simultaneously.  Indeed,
consider the following problem (from, e.g., \cite{Scott06}), which
asks to partition the vertex set of a digraph $D$ into two parts $V_1$
and $V_2$ such that both $e(V_1, V_2)$ and $e(V_2, V_1)$ are large.
It appears that some of the ideas developed in this paper can be used
to attack this problem.  We will pursue this approach further in a
subsequent paper.

\medskip

\noindent \textbf{Acknowledgement.} We would like to thank the anonymous
referees for their valuable comments.


\begin{thebibliography}{100}
\bibitem{AlBoKrSu03}
N.~Alon, B.~Bollob\'as, M.~Krivelevich, and B.~Sudakov,
\newblock Maximum cuts and judicious partitions in graphs without short cycles,
\newblock {\em J. Combin. Theory Ser. B} 88 (2003), 329--346.

\bibitem{AnGrLi83}
L.~Anderson, D.~Grant, and N.~Linial,
\newblock Extremal $k$-colourable subgraphs,
\newblock {\em Ars Combin.} 16 (1983), 259--270.

\bibitem{BoReTh93} B.~Bollob\'as, B. Reed, and A. Thomason
\newblock An extremal function for the achromatic number,
\newblock {\bf Graph Structure Theory}, Eds. N. Robertson and P. Seymour,
\newblock American Mathematical Society, Providence, Rhode Island
(1993), 161--165.

\bibitem{BoSc99}
B.~Bollob\'as and A.~Scott,
\newblock Exact bounds for judicious partitions of graphs,
\newblock {\em Combinatorica} 19 (1999), 473--486.

\bibitem{BoSc00}
B.~Bollob\'as and A.~Scott,
\newblock Judicious partitions of 3-uniform hypergraphs,
\newblock {\em Eur. J. Combinat.} 21 (2000), 289--300.

\bibitem{BoSc02a}
B.~Bollob\'as and A.~Scott,
\newblock Problems and results on judicious partitions,
\newblock {\em Random Structures and Algorithms} 21 (2002), 414--430.

\bibitem{BoSc04}
B.~Bollob\'as and A.~Scott,
\newblock Judicious partitions of bounded-degree graphs,
\newblock {\em J. Graph Theory} 46 (2004), 131--143.

\bibitem{BoSc10}
B.~Bollob\'as and A.~Scott,
\newblock Max $k$-cut and judicious $k$-partitions,
\newblock {\em Discrete Mathematics} 310 (2010), 2126--2139.

\bibitem{Edwards73}
C.~Edwards,
\newblock Some extremal properties of bipartite subgraphs,
\newblock {\em Canad. J. Math.} 25 (1973), 475--485.

\bibitem{Edwards75}
C.~Edwards,
\newblock An improved lower bound for the number of edges in a largest
bipartite subgraph,
\newblock {\em Proc. 2nd Czech Symp. Graph Theory, Prague} (1975), 167--181.

\bibitem{ErGyKo97}
P.~Erd\H{o}s, A.~Gy\'arf\'as, and Y.~Kohayakawa,
\newblock The size of the largest bipartite subgraphs,
\newblock {\em Discrete Math.} 177 (1997), 267--271.

\bibitem{FeLa06}
U.~Feige and M.~Langberg,
\newblock The $RPR^2$ rounding technique for semidefinite programs,
\newblock {\em J. Algorithms} 60 (2006), 1--23.

\bibitem{FrJe97}
A.~Frieze and M.~Jerrum,
\newblock Improved approximation algorithms for MAX $k$-CUT and MAX BISECTION,
\newblock {\em Algorithmica} 18 (1997), 61--77.

\bibitem{GoWi95}
M.~Goemans and D.~Williamson,
\newblock 0.878 approximation algorithms for MAX CUT and MAX 2-SAT,
\newblock {\em Proc. 26th ACM Symp. Theory Comput.} (1994), 422--431;
\newblock Updated as: Improved approximation algorithms for maximum
cut and satisfiability problems using semidefinite programming,
\newblock {\em J. ACM} 42 (1995), 1115--1145.

\bibitem{HaSz70}
A.~Hajnal and E.~Szemer\'edi,
\newblock Proof of a conjecture of P. Erd\H{o}s,
\newblock in {\em Combinatorial theory and its applications, II (Proc. Colloq., Balatonf\"{u}red, 1969)}, North-Holland, Amsterdam (1970), 601--623.

\bibitem{Ha11} J.~Haslegrave, 
  \newblock The Bollob\'as-Thomason conjecture for 3-uniform hypergraphs,
  \newblock {\em Combinatorica} 32 (2012), 451--471.

\bibitem{Hastad01}
J.~H{\aa}stad,
\newblock Some optimal inapproximability results,
\newblock {\em J. ACM} 48 (2001), 798--859.

\bibitem{JLR}
S.~Janson, T.~\L uczak, and A. Ruci\'nski,
\newblock {\bf Random Graphs}, Wiley, New York (2000).

\bibitem{KuOs07}
D.~K\"uhn and D.~Osthus,
\newblock Maximizing several cuts simultaneously,
\newblock {\em Combinatorics, Probability and Computing} 16 (2007), 277--283.

\bibitem{LeTu82}
J.~Lehel and Zs.~Tuza,
\newblock Triangle-free partial graphs and edge-covering theorems,
\newblock {\em Discrete Math.} 30 (1982), 59--63.

\bibitem{Locke82}
S.~Locke,
\newblock Maximum $k$-colorable subgraphs,
\newblock {\em J. Graph Theory} 6 (1982), 123--132.

\bibitem{MaYaYu10}
J.~Ma, P.~Yan, and X.~Yu,
\newblock On several partition problems of Bollob\'as and Scott,
\newblock {\em J. Comb. Theory, Ser. B} 100 (2010), 631--649.

\bibitem{MaYu11}
J.~Ma and X.~Yu,
\newblock Partitioning 3-uniform hypergraphs,
\newblock {\em J. Comb. Theory, Ser. B} 102 (2012), 212--232.

\bibitem{PoTu82}
S.~Poljak and Zs.~Tuza,
\newblock Bipartite subgraphs of triangle-free graphs,
\newblock {\em SIAM J. Discrete Math.} 7 (1994), 307--313.

\bibitem{Scott06}
A.~Scott,
\newblock Judicious partitions and related problems,
\newblock in {\em Surveys in combinatorics}, London Math. Soc. Lecture Note Ser., 327, Cambridge Univ. Press, Cambridge (2005), 95--117.

\bibitem{TSSW00}
L.~Trevisan, G.~Sorkin, M.~Sudan, and D.~Williamson,
\newblock Gadgets, approximation, and linear programming,
\newblock {\em SIAM J. Comput.} 29 (2000), 2074--2097.

\end{thebibliography}
\end{document}